\documentclass[aps, prl, twocolumn, groupedaddress, showpacs, floatfix,10pt]{revtex4-1}

\usepackage{amsmath, amssymb}
\usepackage[english]{babel}
\usepackage[utf8]{inputenc}
\usepackage{xy}
\usepackage{graphicx}
\usepackage{verbatim}
\usepackage{enumerate}
\usepackage{url}
\usepackage{subfigure}
\usepackage{bm}
\usepackage{times}
\usepackage{todonotes}
\usepackage[version=4]{mhchem}
\usepackage{siunitx}

\DeclareSIUnit{\molar}{M}
\DeclareSIUnit{\siemens}{S}

\newcommand{\note}[1]{}

\usepackage{xcolor}
\definecolor{colCB}{rgb}{0,0,.7}
\definecolor{colIK}{rgb}{.7,0,0}
\definecolor{colKW}{rgb}{0,.7,0}
\definecolor{colhl}{rgb}{.9,0,0}

\newcommand{\ee}{\mathrm{e}}

\newcommand{\fig}[1]{Figure~\ref{#1}}
\newcommand{\figpan}[2]{Figure~\ref{#1}{{#2})}}

\graphicspath{{Figures/}}

\newcommand{\maxdim}{N}

\newcommand{\R}{\mathbb{R}}
\newcommand{\Rn}{\R^\maxdim}

\newcommand{\Z}{\mathbb{Z}}

\newcommand{\Tor}{\mathbf{T}}

\begin{document}


\title{Strong coupling yields abrupt synchronization transitions in coupled oscillators}
\author{Jorge L.~Ocampo-Espindola${}^{a}$, Istv\'{a}n Z.~Kiss${}^{a}$, Christian
  Bick${}^{b,d}$, Kyle C.~A.~Wedgwood${}^{b,c}$}%
\email{k.c.a.wedgwood@exeter.ac.uk}

\affiliation{%
${}^a$Department of Chemistry, Saint Louis University,  St.~Louis, MO 63103,
USA \\
${}^b$Centre for Systems Dynamics and Control and Department of Mathematics,
University of Exeter, EX4~4QF, UK \\
${}^c$Living Systems Institute, University of Exeter, EX4~4QD, UK\\
${}^d$Department of Mathematics, Vrije Universiteit Amsterdam, DE Boelelaan 1111, Amsterdam, NL 
}
\date{\today}


\begin{abstract}
\noindent
{
  Coupled oscillator networks often display transitions between
  qualitatively different phase-locked solutions---such as synchrony and
  rotating wave solutions---following perturbation or parameter variation.
  In the limit of weak coupling, these transitions can be understood in
  terms of commonly studied phase approximations.
  As the coupling strength increases, however, predicting the location and
  criticality of transition, whether continuous or discontinuous, from the
  phase dynamics may depend on the order of the phase approximation---or a
  phase description of the network dynamics that neglects amplitudes may
  become impossible altogether.
  Here we analyze synchronization transitions and their criticality
  systematically for varying coupling strength in theory and experiments
  with coupled electrochemical oscillators.
  First, we analyze bifurcations analysis of synchrony and splay states in
  an abstract phase model and discuss conditions under which synchronization 
  transitions with different criticalities are possible.
  Second, we illustrate that transitions with different criticality indeed
  occur in experimental systems.
  Third, we highlight that the amplitude dynamics observed in the
  experiments can be captured in a numerical bifurcation analysis of
  delay-coupled oscillators.
  Our results showcase that reduced order phase models may miss important features that one would expect in the dynamics of the full system.
}
\end{abstract}
\maketitle


\section{Introduction}

Collective oscillatory dynamics are a hallmark of a multitude of real world
networks, such as electrical activity in the brain~\cite{Ashwin2016a,Bick2020},
power grids~\cite{Filatrella2008,Dorfler2013}, and
epidemiology~\cite{Yan2007,Gross2008}. Such systems are often described using
network dynamical systems models that couple together nodes that each
intrinsically (i.e., in the absence of coupling) exhibit stable, hyperbolic
limit cycle oscillations. If the oscillation frequencies of the nodes, or
subsets of nodes, are sufficiently close together, the network can display
phase-locked behaviour in which the phase difference between pairs of nodes
converges to a finite value~\cite{Hoppensteadt97}. As parameters,
such as the coupling strength, are varied, networks may exhibit sharp
transitions between collective oscillations with different phase-difference
properties. Particularly striking examples include the abrupt synchronization
phenomenon in which a group of nodes (potentially encompassing the entire
network) exhibits a sharp transition from an incoherent state to a phase-locked
state in which the phase differences between nodes in the group
vanishes~\cite{Zhang2014,Vlasov2015}.

If the coupling is sufficiently weak, the network dynamics can be described
using phase reduction~\cite{Nakao2015,Pietras2019}. The phase reduction
describes the dynamics of the phases on an attracting invariant torus in terms
of intrinsic rotation and a {phase interaction function} that captures how the
oscillators' phases interact.  To first order, the phase interaction function is
a convolution of the {phase response function}, which captures the linear
sensitivity of the phase of a node oscillation to a perturbation, and a
{coupling function} that describes how nodes interact with one another.  These
functions can often be inferred from data, or estimated using perturbative
experiments. This makes weakly coupled oscillator theory an attractive framework
for studying real world systems, for example to design the dynamics of coupled
oscillator networks~\cite{Kori2008,Kiss2018}.

If the coupling between individual units becomes strong---as is the case in many
real-world systems---the assumptions that underlie phase reduction cease to be
satisfied.  It is thus pertinent to ask which predictions of the weakly coupled
theory break as the coupling strength is increased and how such predictions change.
For example, strong coupling may turn a continuous synchronization transition into a
discontinuous one~\cite{Calugaru2018}.  Recent work has demonstrated that
predictions for infinitesimal coupling strengths are inconsistent with those for
small, finite coupling strengths, even for simple oscillator
models~\cite{Borgers2023}.  Similarly, perturbations to oscillation amplitudes can
impact phase dynamics, particularly, if the amplitudes of different node
oscillations are perturbed in different ways.  For example, it has long been
known that chaotic dynamics with small amplitude variation can emerge close to
bifurcations of coupled oscillator networks as the coupling strength is
increased due to the presence of symmetries in the underlying dynamics in what
is known as {instant chaos}~\cite{Guckenheimer1992}.

To understand network dynamics beyond the weak coupling limit, new mathematical
tools have recently become available.  On the one hand, these include
higher-order phase reductions, that give a more accurate description of the
phase dynamics~\cite{Skardal2020,Bick2023}.  On the other hand, one can derive approximations that
allow for additional degrees of freedom. For example, phase-amplitude reductions
adds a degree of freedom that corresponds to an ``amplitude'' variable;
cf.~\cite{Wedgwood2013,Wilson2016,Letson2018,Wilson2019}.  Such approximations have also
been derived for dynamical systems with time delay~\cite{Kotani2012}.  Despite
being ad-hoc and without a rigorous mathematical justification, there have been
promising results showcasing the merits of these frameworks, including
demonstrations of how amplitude variations can be controlled~\cite{Wilson2020b}.
However, it remains an open question how best to incorporate the effects of
strong coupling in a practical sense.

Here, we take an interdisciplinary approach to elucidate the effect of strong
coupling on the synchronization dynamics in a minimal network of two
delay-coupled phase oscillators.  Specifically, we demonstrate how abrupt
transitions between different phase-locked states of a two-node network are
induced by changes in coupling strength.
First, we consider phase dynamics for two coupled oscillators with higher
harmonics. Here, one would expect higher harmonics to shape the dynamics for highly 
nonsinusoidal oscillations as the coupling strength is increased.
We show that higher harmonics can introduce changes in the criticality of key
bifurcations, which in turn leads to bistability between solutions with
different asymptotic phase differences.  Second, we demonstrate that such
transitions arise in experiments involving a network of
electrochemical oscillators coupled through 
delayed feedback. 
Since the phase theory is insufficient to describe amplitude variations observed in
the experiments, we investigate a network of two delay-coupled oscillators through
numerical bifurcation analysis of a model of the chemical oscillator
network. Here, we demonstrate the existence of branches of symmetry-broken
solutions that are well-matched to the experimental observations.

\section{Continuous and discontinuous transitions between synchronized states in phase oscillators}
\label{sec:Phase}

To understand transitions between in-phase and anti-phase dynamics, we consider the simple case of a network of two delay-coupled nonlinear oscillators. 
Specifically, the state of each oscillator is given by $x_k\in\Rn$ and evolves according to
\begin{subequations}\label{eq:NonlinOsc}
\begin{align}
\dot x_1 &= F(x_1) + K G(x_2(t-\tau)), \\
\dot x_2 &= F(x_2) + K G(x_1(t-\tau)),
\end{align}
\end{subequations}
where~$F:\Rn\to\Rn$ determines the intrinsic oscillatory dynamics and
$G:\Rn\to\Rn$ determines the interactions with strength $K\geq0$ and delay
$\tau\geq 0$. In the uncoupled case, with $K=0$, each node
possesses a stable hyperbolic limit cycle with intrinsic frequency $\omega$.
If the coupling is sufficiently weak ($|K|\ll 1$), the dynamics of
\eqref{eq:NonlinOsc} evolve on an invariant torus in which the oscillator
amplitudes are slaved to the respective oscillator
phases~$\theta_1,\theta_2\in\Tor:=\R/2\pi\Z$.
In this case, the dynamics can be simplified via projection
onto this invariant torus~\cite{Nakao2015,Pietras2019}, so that the (averaged) phase equations for~\eqref{eq:NonlinOsc} with~$n$ relevant harmonics can be written as
\begin{subequations}\label{eq:PhaseOsc}
\begin{align}
\dot \theta_1 &= \omega + g(\theta_2-\theta_1+\alpha),\\
\dot \theta_2 &= \omega + g(\theta_1-\theta_2+\alpha),
\end{align}
\end{subequations}
where
\begin{equation}
g(\phi) = \frac{1}{2}\sum_{m=1}^n{a_j}\sin(m\phi + \gamma_m)
\end{equation}
is the $2\pi$-periodic (phase) coupling function, and~$\alpha$ is a phase shift
parameter common to both oscillators. Up to rescaling, we may assume $a_1 = 1$ and $\gamma_1=0$. 
Note that, in the limit of weak coupling, the delay $\tau$ in \eqref{eq:NonlinOsc} is associated with the phase shift $\alpha$ in \eqref{eq:PhaseOsc}, which, in turn, 
can affect the stability of the synchronized solutions and thus serves as a convenient bifurcation parameter that can be used to engineer 
phase differences between coupled oscillators~\cite{Sakaguchi1986,Kori2008}.

\subsection{Symmetries, bifurcations, and criticality}

By symmetry, the in-phase solution $\Theta_0 =
\{\theta_1=\theta_2\}$ and the anti-phase solution $\Theta_\pi =
\{\theta_1=\theta_2+\pi\}$ are (relative) equilibria of \eqref{eq:PhaseOsc} for any choice of parameter values. 
Note that~\eqref{eq:PhaseOsc} inherits the permutational symmetry $(\theta_1, \theta_2) \mapsto (\theta_2, \theta_1)$ from~\eqref{eq:NonlinOsc}, and---since it describes the slow evolution of the phase differences---a continuous phase-shift symmetry where $\gamma\in\Tor$ acts by $\gamma:(\theta_1, \theta_2)\mapsto (\theta_1+\gamma, \theta_2+\gamma)$.
To eliminate this phase shift symmetry, we can describe the dynamics of~\eqref{eq:PhaseOsc} in terms of the phase difference $\psi:=\theta_2-\theta_1$ between the two oscillators. 
The phase difference evolves according to
\begin{equation}\label{eq:PhaseDiff}
\begin{split}
\dot\psi &= g(-\psi+\alpha)-g(\psi+\alpha)\\
&= \sum_{m=1}^{n}a_m\cos(m\alpha+\gamma_m)\sin(m\psi).
\end{split}
\end{equation}
In-phase synchrony~$\Theta_0$ in~\eqref{eq:PhaseOsc} corresponds to~$\psi=0$ and anti-phase synchrony~$\Theta_\pi$ corresponds to~$\psi=\pi$; both of these points are equilibria of~\eqref{eq:PhaseDiff}.

We now consider bifurcations of in-phase ($\psi=0$) and anti-phase
($\psi=\pi$) configurations as the phase-shift parameter~$\alpha$ is varied.
For coupling functions~$g$ with a single nontrivial harmonic, i.e., $a_j=0$ for
$j\geq 2$, both $\psi=0$ and $\psi=\pi$ bifurcate at
$\alpha=\frac{\pi}{2}+q\pi$, $q\in\Z$ and are connected by a ``vertical'' branch
of equilibria along which any~$\psi\in\Tor$ is an equilibrium of~\eqref{eq:PhaseDiff}. 
If the second harmonic is also non-zero, then there is a nondegenerate branch of equilibria around $\alpha=\frac{\pi}{2}$ that connects $\psi=0$ and $\psi=\pi$~\cite{Rusin2010}; the bifurcations of these solutions are either both super- or both subcritical.
While first and second harmonics may be a suitable approximation in certain
parameter regimes (e.g., where~$G$ describes linear coupling or when the
uncoupled limit cycles are almost sinusoidal in nature), one expects that higher harmonics in the phase dynamics become more relevant in~\eqref{eq:PhaseOsc} as the coupling strength~$K$ is increased. 

It is then instructive to ask what the consequence of the presence of these
higher harmonics might be for the phase dynamics.
One specific important question is whether it is possible for the bifurcations
of the in-phase ($\psi=0$)
and anti-phase configurations ($\psi=\pi$) in~\eqref{eq:PhaseOsc} to have different criticality when higher harmonics are taken into account.
While one can generically control the criticality of the transition locally~\cite{Kuehn2020}, we consider in-phase and anti-phase configurations here simultaneously in the context of~\eqref{eq:PhaseOsc}.

If the harmonics do not have distinct phase shifts, i.e., $\gamma_m = 0$, then the criticality of the bifurcations of $\psi=0$ and $\psi=\pi$ are identical;
this implies in particular that generalizing the phase interaction function considered in~\cite{Rusin2010} to more than two nontrivial harmonics cannot give transitions of distinct criticality.
This can be seen by noting that the system for $\gamma_m = 0$ has a parameter symmetry $(\psi,\alpha)\mapsto(\pi-\psi,\pi-\alpha)$.
This implies that $\psi=\frac{\pi}{2}$ is an equilibrium for $\alpha = \frac{\pi}{2}$ and that any bifurcation of $\psi=0$ at~$\alpha = \hat\alpha$ leads to an identical bifurcation of $\psi=\pi$ at~$\alpha = \pi-\hat\alpha$. 
Moreover, if all even harmonics vanish (i.e., $a_m$ = 0 for $m$ even), then $\psi=0$ and $\psi=\pi$ bifurcate at $\alpha = \frac{\pi}{2}$.
Thus, if parameters are such that there is only a single bifurcation of
$\psi=0,\pi$ for $\alpha\in(0,\pi)$ (i.e., these bifurcations are related by
symmetry) then it is necessary to have non-zero~$\gamma_m$ for the bifurcations to have distinct criticality.

\subsection{Distinct criticality of transitions of in-phase and anti-phase configurations}

We now consider the bifurcations of $\psi\in\{0,\pi\}$ in~$\alpha$ for more general choices of~$\gamma_m$. 
Expanding~\eqref{eq:PhaseDiff} around $\psi=0$ yields
\begin{align}
\dot\psi = \sum_{m=1}^{n}a_m\cos(m\alpha+\gamma_m)\left(m\psi - \frac{m^3\psi^3}{3!} + \dotsb\right).
\end{align}
Thus, the linear stability of $\psi=0$ as well as the criticality of the
(potentially degenerate) pitchfork bifurcation around $\psi=0$ are determined by
\begin{subequations}\label{eq:SyncBif}
\begin{align}
D_1(0; \alpha) &= \sum_{m=1}^{n}ma_m\cos(m\alpha+\gamma_m)\label{eq:BifZero}, \\
D_3(0; \alpha) &=
-\frac{1}{3!}\sum_{m=1}^{n}m^3a_m\cos(m\alpha+\gamma_m)\label{eq:CritZero}.
\end{align}
\end{subequations}
Specifically, solving $D_1(0; \alpha)=0$ for~$\alpha$ determines a bifurcation
point~$\alpha^{(0)}$ of $\psi=0$ as~$\alpha$ is varied and $D_3(0;
\alpha^{(0)})$~determines the criticality of this transition. In particular, 
the bifurcation yields a continuous transition (a supercritical pitchfork bifurcation with an emerging branch of stable equilibria) if $D_3(0; \alpha^{(0)})<0$ and a discontinuous transition (a subcritical pitchfork bifurcation with an emerging branch of unstable equilibria) if $D_3(0; \alpha^{(0)})>0$.
In a similar fashion, expanding~\eqref{eq:PhaseDiff} around $\psi=\pi$ gives
\begin{subequations}\label{eq:SplayBif}
\begin{align}
D_1(\pi; \alpha) &= \sum_{m=1}^{n}(-1)^{m} ma_m\cos(m\alpha+\gamma_m), \\
D_3(\pi; \alpha) &=
-\frac{1}{3!}\sum_{m=1}^{n}(-1)^{m}m^3a_m\cos(m\alpha+\gamma_m)\label{eq:CritPi}.
\end{align}
\end{subequations}
Thus, the bifurcation of the anti-phase configuration $\psi=\pi$ at~$\alpha^{(\pi)}$ is continuous if $D_3(\pi; \alpha^{(\pi)})<0$ and discontinuous if $D_3(\pi; \alpha^{(\pi)})>0$.

\begin{figure*}
\includegraphics[width=\linewidth]{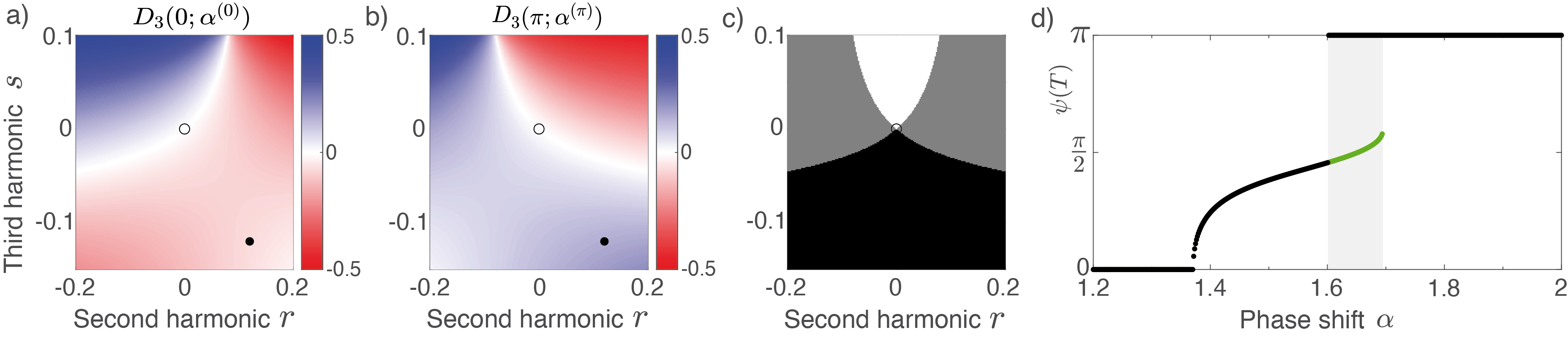}
\caption{\label{fig:Crit}
Varying parameters demonstrates regions in which the bifurcations of equilibria
$\psi=0,\pi$ have distinct criticality.
Results are demonstrated for fixed $\gamma_2 = 0.2$, $\gamma_3 = 0.5$.
a) Criticality coefficient~\eqref{eq:CritZero} for $\psi=0$ as the strength~$r,s$ of the second and third harmonic is varied; blue indicates a continuous, red a discontinuous transition; a hollow circle indicates $r=s=0$ and a filled circle the parameter values in panel d).
b) Criticality coefficient~\eqref{eq:CritPi} for $\psi=\pi$ as in panel a).
c) Regions in which the bifurcations of the in-phase and anti-phase
configurations have distinct criticalities: white indicates that the transition at $\psi=0$ is continuous while $\psi=\pi$ is discontinuous and black vice versa. Bifurcations of $\psi=0$ and $\psi=\pi$ have the same criticality in the grey regions.
d) Pseudocontinuation plot for $r=-s=0.12$ as the
parameter~$\alpha$ is increased (green) or decreased (black); an approximate
region of hysteresis/multistability due to the discontinuous transition of the
$\psi=\pi$ solution is shaded in grey.
}
\end{figure*}

We now show explicitly that there is an open set of parameters for which the criticalities of the transition of the in-phase and anti-phase configurations are distinct.
We first restrict to coupling functions whose first three harmonics may be nontrivial ($a_2=r$, $a_3=s$, $a_m=0$ for $m>3$).
Recall that for $r=s=0$, the equilibria $\psi=0$, $\psi=\pi$ undergo a
degenerate bifurcation at $\alpha=\frac{\pi}{2}$ with a ``vertical'' bifurcation
branch (i.e., any~$\psi\in\Tor$ is an equilibrium). 
For~$r,s$ small, this branch will be perturbed which leads (generically) to nondegenerate pitchfork bifurcations at $\alpha\approx\frac{\pi}{2}$.
Let $\beta=\frac{\pi}{2}-\alpha$ denote the deviation of the bifurcation point from $\alpha=\frac{\pi}{2}$.
Assuming that $m\beta+\gamma_m$ is small, we can approximate these bifurcation
points by
expanding the cosine term in~\eqref{eq:BifZero} and collecting terms in~$\beta$
up to first order to give the approximate location of the  bifurcation point of $\psi=0$ as
\begin{align}
\tilde \alpha^{(0)} &=  
\frac{\pi}{2} - \frac{a_2+a_1\gamma_1-a_3\gamma_3}{a_1-3a_3} = \frac{\pi}{2}-\frac{r-s\gamma_3}{1-3s}.
\intertext{Using the same approximation with $\tilde \beta^{(\psi)}=\frac{\pi}{2}-\tilde \alpha^{(\psi)}$, the criticality at the approximate bifurcation point is}
\begin{split}
\widetilde C^{(0)} &=  a_1\gamma_1+3^2a_2-3^3a_3\gamma_3+(3^4a_3-a_1)\tilde \beta^{(0)}\\ 
&= 9r-27s\gamma_3+(81s-1)\tilde \beta^{(0)}.
\end{split}
\intertext
{Similarly, we can approximate the bifurcation point of $\psi=\pi$ by}
\tilde \alpha^{(\pi)} &=  \frac{\pi}{2}-\frac{r+s\gamma_3}{3s-1}
\intertext{with criticality determined by}
\widetilde C^{(\pi)} &=  9r+27s\gamma_3+(1-81s)\tilde \beta^{(\pi)}.
\end{align}

To see that there is an open set of parameters for which the criticality of
$\psi=0$ and $\psi=\pi$ is distinct, consider the case with a vanishing second
harmonic, $r=0$. Then, $\beta^{(0)}=\beta^{(\pi)}$ and 
\[\widetilde C^{(0)} = \left(\frac{81s-1}{3s-1}-27\right)\gamma_3s = -\widetilde C^{(\pi)}.\]
Thus, for $\widetilde C^{(0)}\neq 0$, the transitions of in-phase and anti-phase
configurations have distinct criticality, which are exchanged as~$s$ passes through zero.
Since the expressions considered are continuous in all parameters for small~$s$,
this yields an open set of parameters for which the in-phase and anti-phase configurations have
distinct criticality, as claimed above. Note that this phenomenon is not limited
to the case with three harmonics with parameters~$r,s$ but also occurs if we allow small non-zero~$a_m$, $m>3$.

To demonstrate our findings, we compute the bifurcation points and their criticality numerically using~\eqref{eq:SyncBif} and~\eqref{eq:SplayBif};  
cf.~\figpan{fig:Crit}{a,b}.
There is indeed an open set of parameters for which the bifurcations of the
in-phase and anti-phase configurations have distinct criticality, as shown in \figpan{fig:Crit}{c}. 
For a slowly varying parameter~$\alpha$~\footnote{For each
  parameter~$\alpha_n$, we solved~\eqref{eq:PhaseDiff}, numerically for $T=5,000$
  time units with initial condition being~$\psi(T)$ for parameter $\alpha_{n\pm
  1}$ plus a small random perturbation.}, this results in the bifurcation behavior shown in \figpan{fig:Crit}{d} where the transition of $\psi=0$ at the bifurcation point is continuous while $\psi=\pi$ shows a discontinuous, abrupt transition. 
Note that we here focus on the bifurcation points that converge to $\alpha=\frac{\pi}{2}$ as $r,s\to 0$; further bifurcations---also along the nontrivial branch---can occur as the influence of second and third harmonic grows.

\section{Hysteresis in coupled electrochemical reactions}
\label{sec:Exp}

We next investigate the consequence of the results of the preceding section in a
real-world system.  In particular, we examine whether the regions of existence of pitchfork
bifurcations with distinct criticality predicted in \fig{fig:Crit} can be
induced in a network of two oscillatory electrochemical reactions coupled with
time-delayed linear feedback.  Here, we predict that increasing the coupling
strength between the reactions can drive changes in pitchfork criticality and
hence give rise to bistability between phase-locked solutions with different
phase differences.

\subsection{Experimental setup}
\begin{figure}
  \includegraphics[width=0.95 \linewidth]{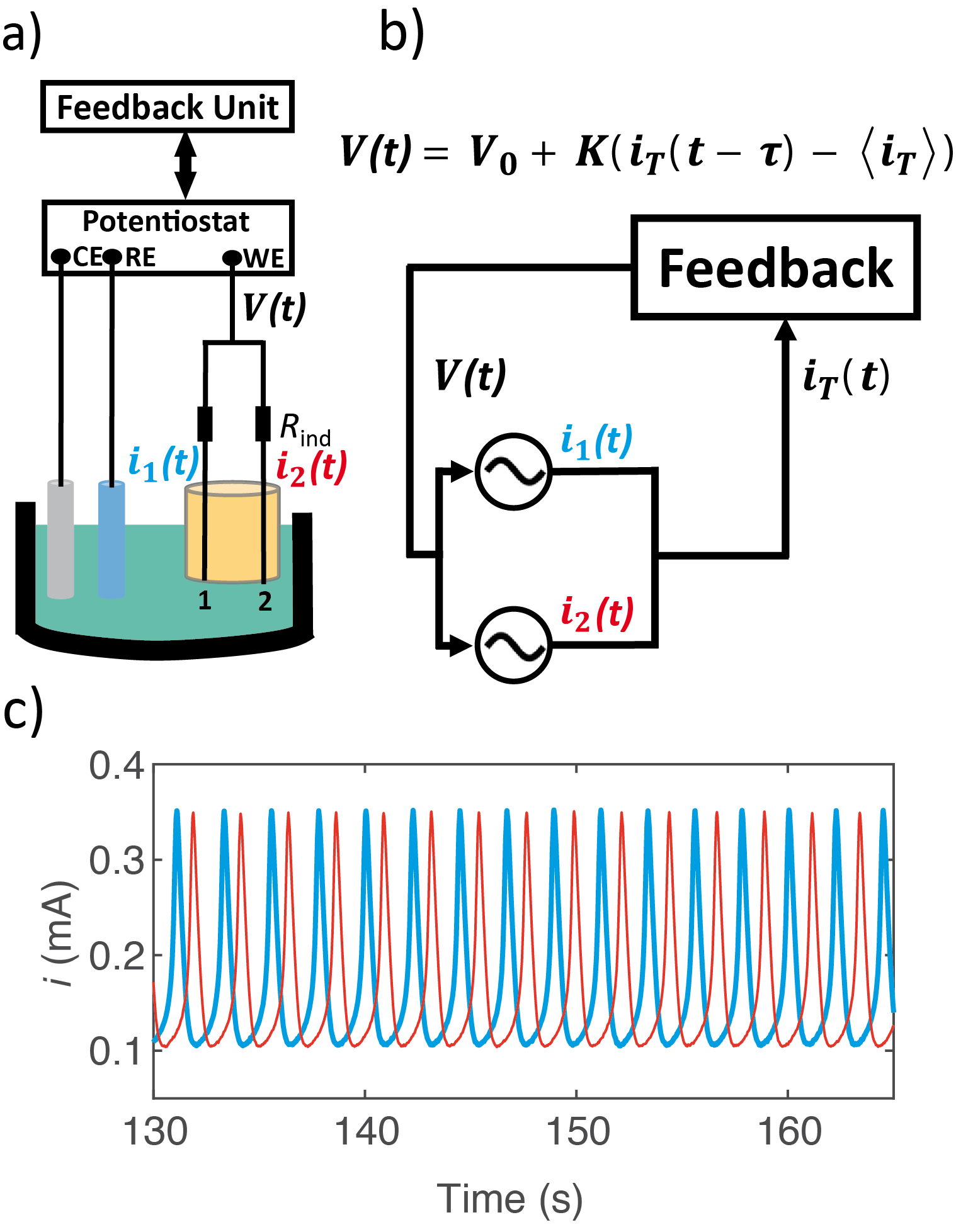}
\caption{\label{fig:exp_set} Illustration of the time-delayed linear feedback
experiment with a time series of the uncoupled system. a) Schematic of the
experimental setup. CE: Counter electrode, RE: Reference electrode, and WE: Working
electrodes. b) Diagram of the delay feedback schematic of the electrochemical
experiment. The currents ($i_1$, $i_2$) of each nickel wires were measured and
added to obtain a total current ($i_T$). The $i_T$ was fed back with a coupling
strength ($K$), a delay, ($\tau$) and applied to the circuit potential ($V(t)$).
c) Time series of the currents for the uncoupled ($K=0$) oscillators and
without delay ($\tau = 0$). The blue and red lines correspond to oscillator 1 and 2, respectively.}
\end{figure}

A schematic of the experimental setup is shown in \figpan{fig:exp_set}{a}.  The
three-electrode electrochemical cell is equipped with a \ce{Pt} coated \ce{Ti}
rod as a counter (C), a \ce{Hg/Hg2SO4/}sat.~\ce{K2SO4} as a reference (R) and
two \ce{Ni} wires (Goodfellow Cambridge Ltd, 99.98\%, 1.0 mm diameter) as
working electrodes (W) connected to a potentiostat (ACM Instruments, Gill AC).
The electrodes are immersed in a \SI{3}{\molar} \ce{H2SO4} solution as an
electrolyte and kept at a constant temperature of \SI{10}{\degreeCelsius}. 
 

When a constant circuit potential with respect to the reference electrode
($V_0=\SI{1200}{\milli\volt}$) is applied by the potentiostat and external
resistance ($R_\mathrm{ind}=\SI{1}{\kilo\ohm}$) is attached to each nickel wire,
the electrochemical dissolution of nickel exhibited periodic oscillations of the
current \cite{kiss_tracking_2006} (see \figpan{fig:exp_set}{c}). 
In our specific experiments, the natural (uncoupled) frequencies of oscillators
1 and 2 were $\omega_1 = \SI{0.446}{\hertz}$ and $\omega_2 =
\SI{0.444}{\hertz}$, respectively, with a mean frequency of \SI{0.445}{\hertz}
and a mean period of $T =\SI{2.25}{\second}$. 


The potentiostat is interfaced with a real-time LabVIEW controller, and is
used to measure the total current $i_T$ and subsequently set the circuit potential $V(t)$ at a rate of 200 Hz according to the equation
\begin{equation}\label{eq:exp_feedback}
V(t) = V_0 + K(i_T(t-\tau)- \langle i_T\rangle),
\end{equation}
where $V(t)$ and $V_0$ are the applied and the offset
circuit potential,
respectively, $K$ is the coupling strength, $i_T$ is the  time averaged total current, $\langle
i_T \rangle$ is the mean value of the total current, and $\tau$ is the time
delay. The coupling between electrodes is induced using external global
feedback via a small adjustment of the circuit potential according to the scheme
in \figpan{fig:exp_set}{b}.
In contrast to previous studies in which nonlinear feedback was used to couple
oscillators very close to a Hopf bifurcation \cite{Rusin2010,
bick2017robust}, the oscillators under consideration here are far from the
Hopf bifurcation point, and are coupled through linear feedback. In particular,
the uncoupled oscillator waveforms are far from the single harmonic profiles
expected for systems close to a Hopf bifurcation. As a result, we predict that
higher harmonics will be important in determining the phase dynamics as the
coupling strength $K$ is increased, in line with results illustrated in
\fig{fig:Crit}.

\subsection{Results with weak coupling}

\begin{figure}
  \includegraphics[width=0.95 \linewidth]{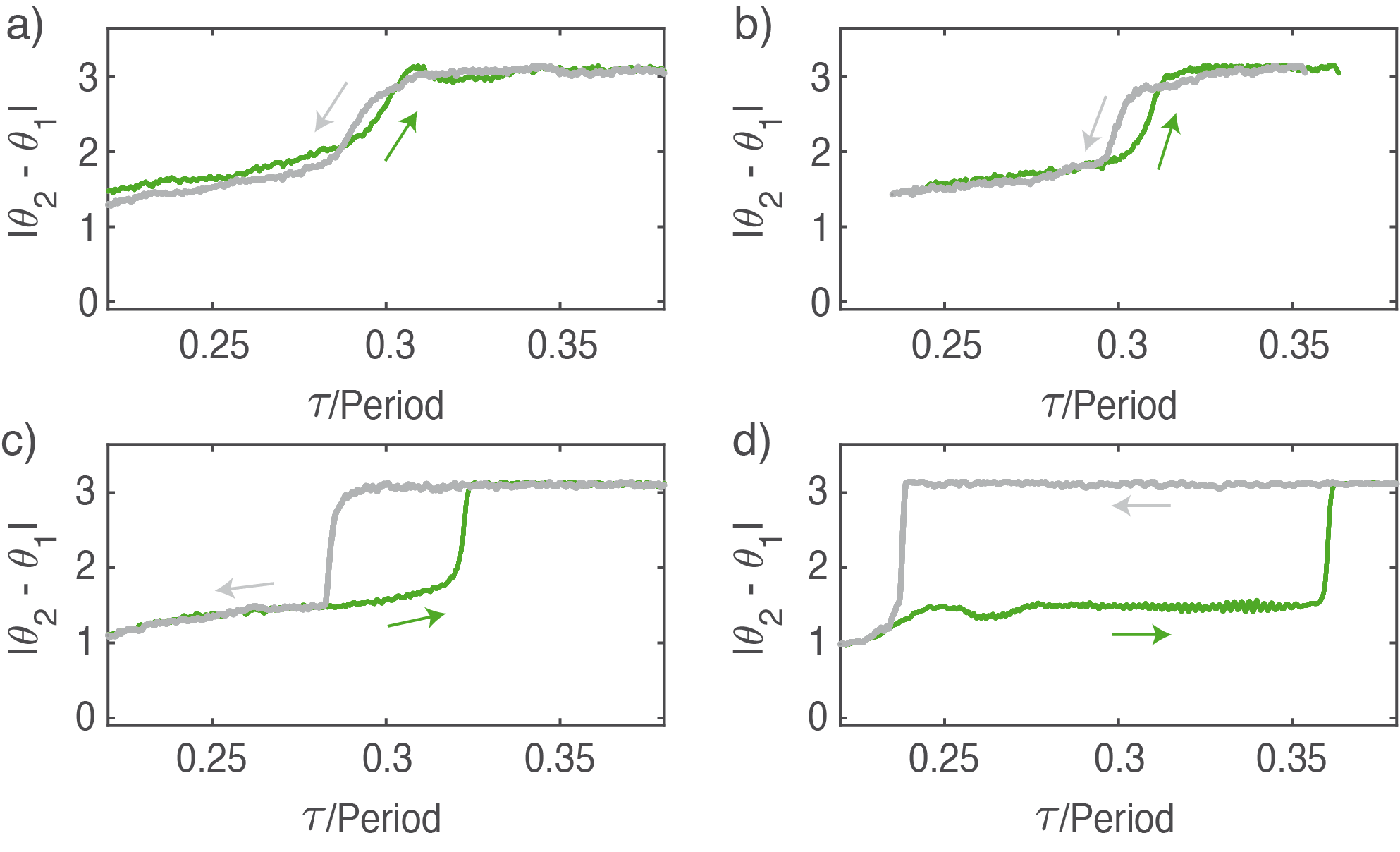}
\caption{\label{fig:slow_scan} Scan under variation of the time delay $\tau$ at
weak coupling strength ($K=\SI{-0.12}{\volt\per\milli\ampere}$) of the two
electrode system. a) Time series of the current for the in-phase behavior for
$\tau=\SI{0.081}{\second}$ (0.036 $\tau/T$). b) Time series of the
current for the out-of-phase behavior for $\tau=\SI{0.60}{\second}$ (0.27
$\tau/T$). c) Time series of the current for the anti-phase behavior for
$\tau=\SI{0.89}{\second}$ (0.39 $\tau/T$). In panels a)-c), the blue
and red lines correspond to oscillators 1 and 2, respectively. d) Phase
difference of the coupled oscillators as a function of the time delay. The green
line is the phase difference corresponds to the forward scan ($\tau$ = 0 $\to$ 0.5
$\tau/T$), and the grey line to the backward scan ($\tau$ = 0.5 $\tau/T$ $\to$
0) with the direction indicated by the green and grey arrows, respectively.}
\end{figure}

We first demonstrate the system dynamics when the coupling is weak
($K=\SI{-0.12}{\volt\per\milli\ampere}$) for different time-delay values. For
illustration, we chose three different time-delays $\tau=\SI{0.081}{\second}$,
$\tau=\SI{0.60}{\second}$, and $\tau=\SI{0.89}{\second}$.  When
$\tau=\SI{0.081}{\second}$, the current signal for each oscillator overlap,
yielding an in-phase synchronized configuration with nearly identical peak-to-peak
amplitudes  ($\Delta A = A_2 - A_1=\SI{-1.0E-3}{\milli\ampere}$, see
\figpan{fig:slow_scan}{a}).  As shown in \figpan{fig:slow_scan}{b}, when the time delay
is increased to $\tau=\SI{0.60}{\second}$, an out-of-phase synchronized configuration
($\Delta \phi = 4.34$ rad) is observed with a relatively large amplitude
difference ($\Delta A = \SI{-0.011}{\milli\ampere}$).  \figpan{fig:slow_scan}{c}
shows the dynamics when we further increased the delay to
$\tau=\SI{0.89}{\second}$ where we observe that the elements synchronized in an
anti-phase configuration with both oscillators having similar amplitudes ($\Delta
A = \SI{-7.0E-4}{\milli\ampere}$).


The quasi-stationary phase difference between the two coupled oscillators
was experimentally measured under slow variation of $\tau$. After letting the
oscillators settle to a synchronized configuration for $\tau=\SI{0.00}{\second}$,
measurements were taken as the time delay was slowly increased to $\tau=\SI{1.12}{\second}$ ($\tau \approx T/2$;
around one half period).
Following this, the time delay was decreased from $\tau=\SI{1.12}{\second}$ back
down to $\tau=\SI{0.00}{\second}$ at the same rate as the forward (increasing
$\tau$) scan.

The phase difference for weak coupling strength,
$K=\SI{-0.12}{\volt\per\milli\ampere}$, is shown in \figpan{fig:slow_scan}{d}. For
time delays $\tau \le 0.1T$, the oscillators exhibit a phase difference close
to  0 (equivalently $2\pi$), indicating in-phase synchronization.  For $0.1T  <
\tau \le 0.3T$, the phase difference between the oscillators increases
monotonically with respect to $\tau$ until it reached a phase difference close
to $\pi$. The oscillators remain anti-phase synchronized when the time delay is
further increased ($0.3T < \tau \le 0.5T$).  For decreasing $\tau$ from
$\tau=0.5T$ to $\tau=0$, the system passes through a sequence of phase-locked
configurations, from anti-phase, to out-of-phase, and finally to in-phase dynamics.  The
green and the grey lines in \figpan{fig:slow_scan}{d} correspond to the scan where
time delay was increased and decreased, respectively, and it can be seen that
the curves approximately overlap. In other words, the transition from in-phase
to anti-phase through out-of-phase synchronized configurations occurs without
hysteresis. 

\subsection{Results with strong coupling}

\begin{figure}
  \includegraphics[width=0.95 \linewidth]{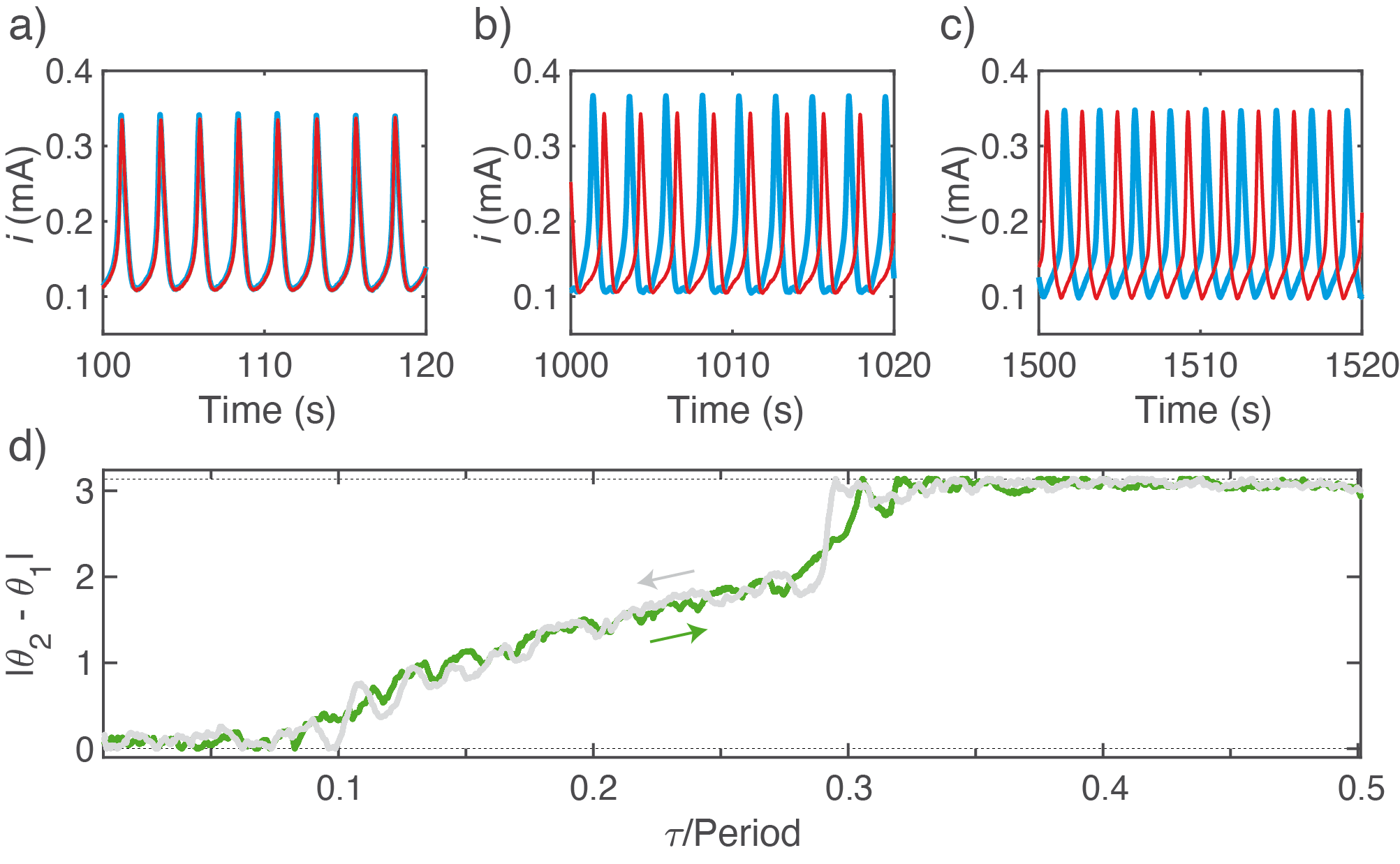}
\caption{\label{fig:zoom_scan} Phase difference under slow variation of the time
delay close to the anti-phase solution for increasing coupling strengths. a)
$K=\SI{-0.12}{\volt\per\milli\ampere}$, b)
$K=\SI{-0.18}{\volt\per\milli\ampere}$, c)
$K=\SI{-0.25}{\volt\per\milli\ampere}$ and d)
$K=\SI{-0.50}{\volt\per\milli\ampere}$. The green line corresponds to the phase
difference in the forward (green arrow) scan and the grey line to the backward
(grey arrow) scan.}
\end{figure}

We next investigated how the phase differences changed as~$\tau$ was increased
and decreased for different coupling strengths, as reported in
\fig{fig:zoom_scan}. When the coupling strength is weak
($K=\SI{-0.12}{\volt\per\milli\ampere}$), the curves for increasing~$\tau$ and
for decreasing $\tau$ overlap, as observed in \figpan{fig:zoom_scan}{a}. 

An increase in the coupling strength to $K=\SI{-0.18}{\volt\per\milli\ampere}$
(see \figpan{fig:zoom_scan}{b}) reveals a small region around $0.30T<\tau <0.31$
where the system possesses two stable stationary configurations coexisting
simultaneously. In this case, the system exhibits bistability, and the curves
for increasing and decreasing $\tau$ do not overlap.  When the coupling strength
increases further to $K=\SI{-0.25}{\volt\per\milli\ampere}$ (see
\figpan{fig:zoom_scan}{c}), the bistability region enlarges to $0.28T<\tau <0.32T$.
Finally, for $K=\SI{-0.50}{\volt\per\milli\ampere}$ (see \figpan{fig:zoom_scan}{d}),
the bistability region is larger still ($0.24T<\tau<0.36T$), resulting in an
extended and well-defined region where the out-of-phase and anti-phase synchronized
configurations co-exist.

\begin{figure}
  \includegraphics[width=0.95 \linewidth]{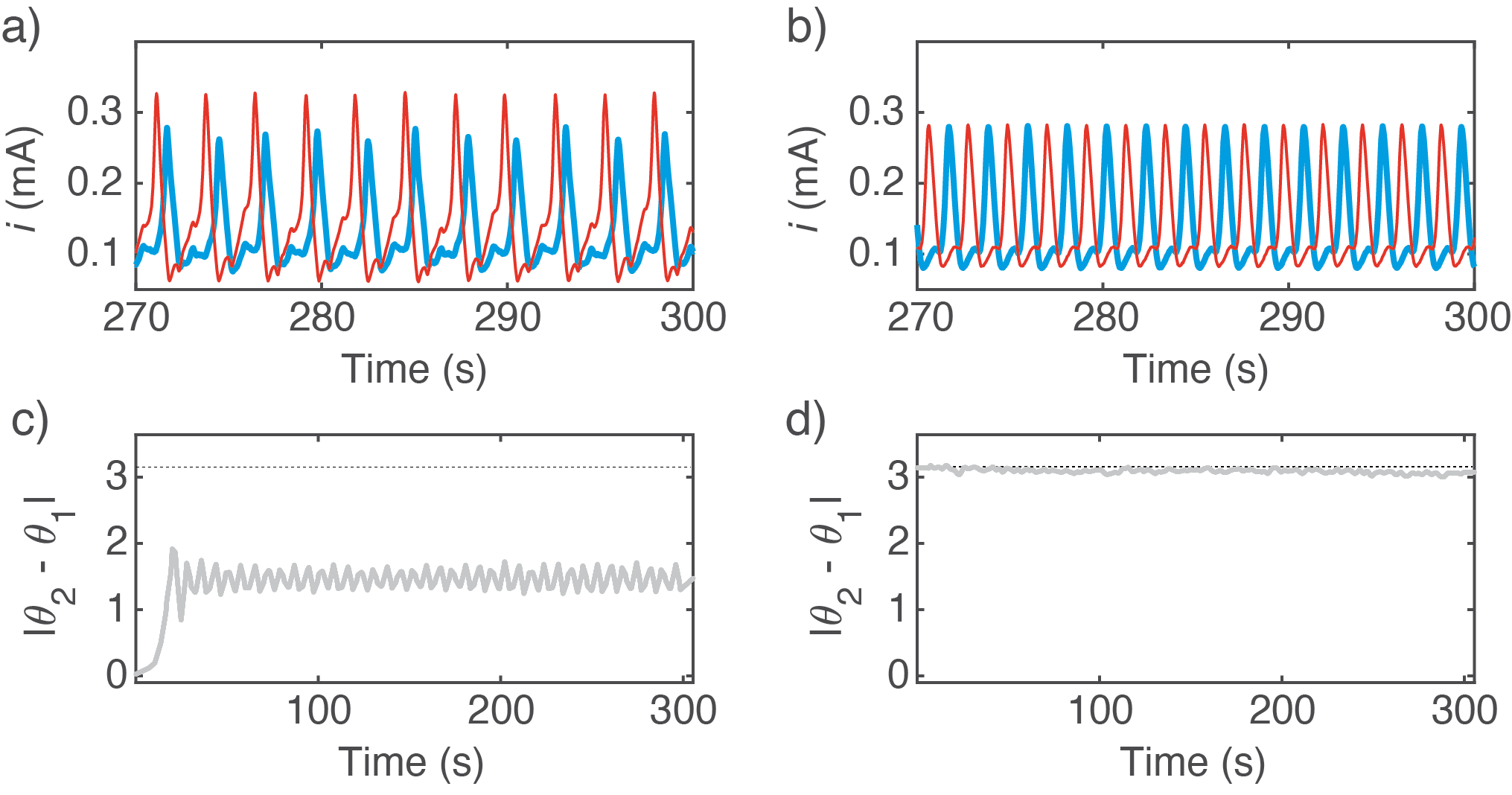}
\caption{\label{fig:phase_diff} Time series and phase difference of the currents
in the bistability region for strong feedback,
$K=\SI{-0.50}{\volt\per\milli\ampere}$, with time-delay
$\tau=\SI{0.70}{\second}$ (0.31 $\tau/T$). The panels a) and c) corresponds to
the out-of-phase configuration, while panels b) and d) correspond to the anti-phase
configuration. The blue and red lines correspond to oscillator 1 and 2, respectively.}
\end{figure}

To better exemplify the bistable nature of the stationary configurations, we performed
experiments in which the system exhibited the bistability phenomena at a strong
feedback gain value $K=\SI{-0.50}{\volt\per\milli\ampere}$ with appropriate initial
conditions (in-phase or anti-phase) and time delay $\tau=\SI{0.70}{\second}$
($\tau=0.31T$).  The time series of the current and the phase difference are
shown in \figpan{fig:phase_diff}{a,b} for an experiment in which the system was
initiated from an in-phase initial condition. After a transient time of about
\SI{25}{s}, the two oscillators transition to an out-of-phase synchronized configuration
with  $\Delta \phi =$ 4.82 rad. Similar to the previous examples, the
out-of-phase synchronized configuration has a relatively large amplitude difference, in
this case, $\Delta A$ = \SI{0.034}{\milli\ampere}.  The corresponding
experimental results starting from anti-phase initial conditions are shown in
\figpan{fig:phase_diff}{c,d}.  As expected, the system remains in the
anti-phase synchronized configuration with a very small amplitude difference ($\Delta A
=$ \SI{1E-4}{\milli\ampere}).

We thus see that electrochemical oscillators display both out-of-phase and
anti-phase configuration for a strong value of the coupling strength
($K=\SI{-0.50}{\volt\per\milli\ampere}$) for different initial conditions,
further confirming the bistability phenomena observed in the bifurcation diagram
in \figpan{fig:zoom_scan}{d}.  The experiments in \fig{fig:phase_diff} also
demonstrate that these configurations remained stable for at least \SI{300}{s} (133 cycles).
We next investigate whether such phenomenon can be attributed to differences in
the criticality of the bifurcations of the in-phase and the anti-phase configurations.

\section{Amplitude asymmetry in a coupled nonlinear oscillator model}
\label{sec:Bif}

The analysis of the phase model \eqref{eq:PhaseDiff} in \fig{fig:Crit} predicts
regions in parameter space in which the pitchfork bifurcations of the
in-phase and anti-phase synchronized solutions have different criticalities. In these
regions, we would expect the bistability between one of these solution types and
an out-of-phase solution, as observed in \fig{fig:zoom_scan}. However, since the
phase model disregards information about oscillation amplitude, it cannot
predict the amplitude asymmetry observed in \figpan{fig:slow_scan}{b} and
\figpan{fig:phase_diff}{a}.
Our goal in this section is to explore the qualitative asymptotic phase dynamics
expected in the electrochemical experiments via bifurcation analysis of a suitable
system of DDEs to further investigate this amplitude asymmetry.
Some of salient synchronization features of the two electrode system have been shown to be well captured by the network Brusselator
model~\cite{Rusin2010}:
\begin{subequations}
  \label{eq:bruss}
  \begin{align}
    \dot{x}_i &= (B-1)x_i + A^2 x_i + f(x_i,y_i) + K G(x_j), \\
    \dot{y}_i &= -Bx_i - A^2 y_i - f(x_i,y_i),
  \end{align}
\end{subequations}
for $j \neq i$, where $i = 1,2$ and
\begin{equation}
  f(x,y) = (B/A)x^2 + 2Axy + x^2y.
\end{equation}
We identify the $x$ component of \eqref{eq:bruss} with the
currents measured in the potentiostat experiments and $y$ with an unobserved
recovery variable. The parameters dictating the intrinsic oscillator dynamics are
hereon set to $A = 0.9$ and $B = 2.3$. For these parameter values and with the
global coupling strength set to 0, each oscillator
possesses a stable hyperbolic limit cycle with period $T =
7.33$. The coupling function, which applies only to the $x$ equations of
\eqref{eq:bruss}, is given by
\begin{equation}
  G(x) = \sum_{n=1}^{2} k_n x(t-\tau_n-\tau)^n.
\end{equation}
We set the amplitude ($k_n$) and delay ($\tau_n$) parameters using the synchronization engineering methods
outlined in \cite{Kori2008}. Briefly, we express the phase response curve of the
uncoupled oscillators as a Fourier series $Z(\theta) = \sum_{n\in\mathbb{N}} Z_n
\mathrm{e}^{in\theta}$ and a target phase interaction function as $g(\theta) =
\sum_{n\in\mathbb{N}} g_n \mathrm{e}^{in\theta}$. The $k_n$ and $\tau_n$
parameters are then chosen so that the Fourier series representation of the coupling
function, i.e., $G(\theta) = \sum_{n\in\mathbb{N}} G_n \mathrm{e}^{in\theta}$,
are approximated by $G_{n} = Z_{n} g_{-n}$. In this study, we use a
phase-shifted Hansel--Mato--Meunier type interaction function given by~\cite{Hansel1993}
\begin{equation}
  \begin{split}
    g(\theta) &= \sin(\theta-\tau) - r \sin(2(\theta -
    \tau)) \\
              &= -\frac{i}{2} \ee^{i\tau} \ee^{-i\theta} + \frac{ir}{2}
              \ee^{2i\tau}\ee^{-2i\theta} + \text{cc.},
  \end{split}
  \label{eq:HMM}
\end{equation}
where~$r$ scales the
contribution of the second harmonic and~$\tau$ is a common phase shift
parameter. We set $r = 0.5$ and consider the system dynamics under variation of
$\tau$ and $K$.

For small $K$, the system dynamics is well approximated by a phase reduced
model of the type given by \eqref{eq:PhaseOsc} and so
the phase difference $\psi$ between the two oscillators obeys \eqref{eq:PhaseDiff}.
Since our choice for $g$ contains only two harmonics, we would not expect the
pitchfork bifurcations of the in-phase and anti-phase solutions to have
different criticalities for small $K$, in contrast to the predictions for
the phase interaction function with three harmonics in \fig{fig:Crit}. In fact,
it has previously been shown experimentally that phase synchronization patterns matching
those expected via \eqref{eq:PhaseDiff} with \eqref{eq:HMM} can be achieved in
the electrochemical experiment. In particular, phase locking with arbitrary steady state phase
differences can be realised through variation of the common delay~$\tau$~\cite{Rusin2010}.
Moreover, the pitchfork bifurcations of the in-phase and anti-phase synchronized solutions are
both supercritical in nature and hence the system does not exhibit any
bistability, unlike that observed in the experiments in \fig{fig:zoom_scan}.

We use the Matlab-based package DDE-BIFTOOL to explore the asymptotic system
dynamics under variation of~$\tau$ as~$K$ is increased. DDE-BIFTOOL is designed
to perform numerical bifurcation and stability analysis of systems with fixed
discrete and/or state-dependent delays. It allows for flexible encoding of
systems and for the specification of additional system constraints, such as
relationships between delays, which we shall leverage to implement the common
delay term. The pipeline for numerical bifurcation analysis is outlined in the
Appendix.

\subsection{Numerical bifurcation analysis results}

\begin{figure}
  \includegraphics[width=\linewidth]{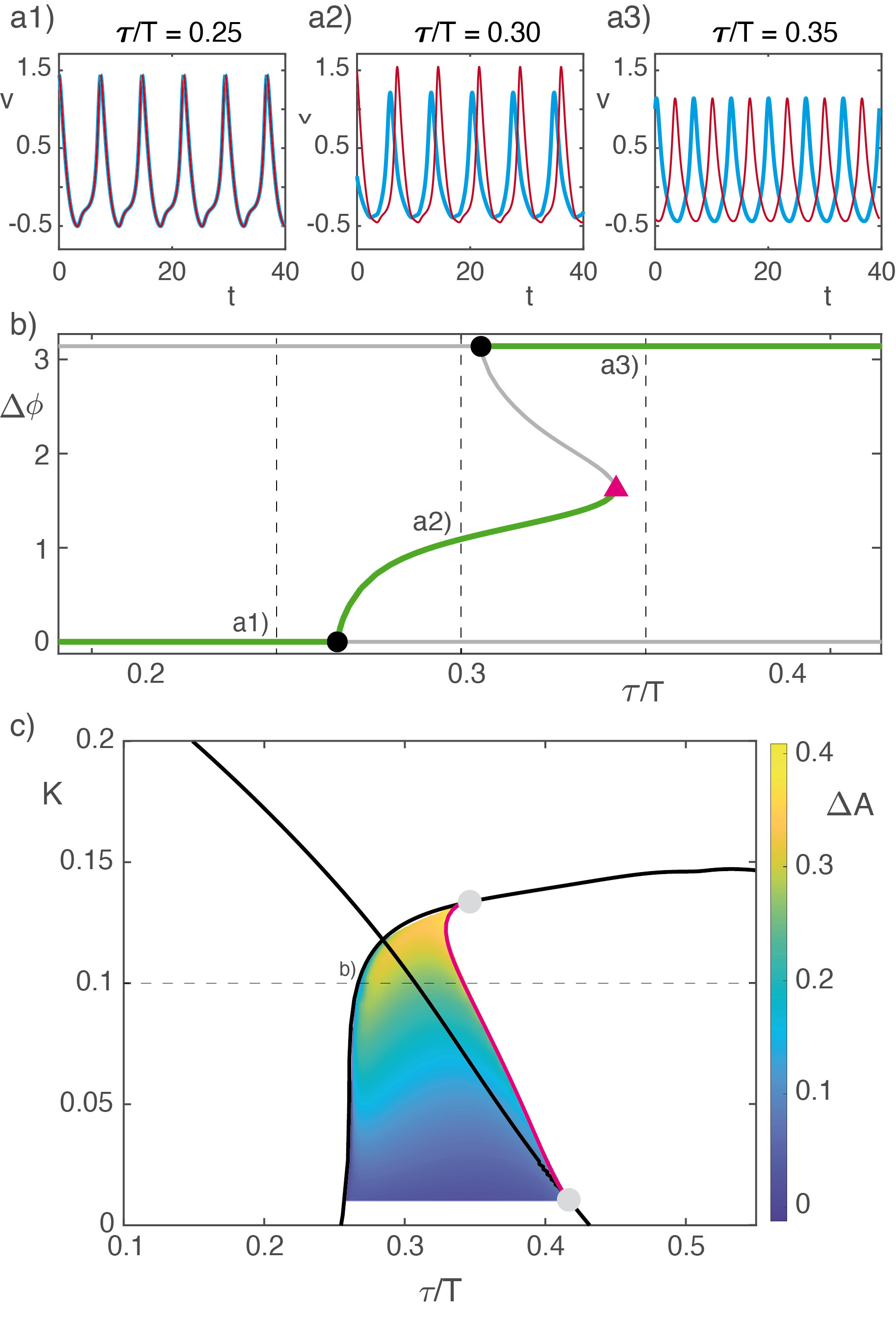}
  \caption{\label{fig:bruss} Bifurcation diagram and time series for the coupled
  Brusselator system \eqref{eq:bruss} computed using DDE-BIFTOOL.
  a1)-a3) Time series solutions corresponding to the in-phase (a1),
  out-of-phase (a2) and anti-phase (a3) synchronized configurations. Blue and red lines
  correspond to oscillator 1 and 2, respectively. 
  b) Bifurcation diagram under variation of $\tau$. Green and grey curves
  correspond to stable and unstable periodic solutions, respectively. Pitchfork bifurcations are
  depicted with black circles, whilst the fold of periodic orbits is marked with a
  pink triangle. The dashed vertical lines show where the corresponding time
  series in panel a) were computed.
  c) Two parameter bifurcation
  diagram under simultaneous variation of $K$ and $\tau$. The black
  curves correspond to pitchfork bifurcations and the pink curve represents the
  fold of periodic orbits. Superimposed on the diagram is a heat map showing the
  relative amplitude asymmetry between the oscillators along the stable part of
  the out-of-phase branch. The dashed horizontal line indicates
  where the bifurcation diagram in panel b) was computed. Bifurcations associated with changes of
  criticality of the pitchforks are depicted by grey markers.}
\end{figure}

The results of the bifurcation analysis procedure are shown in \fig{fig:bruss}. Specifically, \figpan{fig:bruss}{a} showcases temporal profiles of the~$x_i$ along the one parameter
bifurcation diagram shown in \figpan{fig:bruss}{b}. These panels are to be compared with the
equivalent panels in \fig{fig:slow_scan} and \fig{fig:zoom_scan}. \figpan{fig:bruss}{b}
highlights the presence of unstable out-of-phase synchronized solutions along the central
branch. This unstable portion of branch is generated following a change in
criticality of the pitchfork bifurcation of the anti-phase solution. This can be
seen more clearly in the two parameter bifurcation diagram in \figpan{fig:bruss}{c}, where
we observe that the pitchfork of the anti-phase solution becomes subcritical at a
small positive value of~$K$.  This panel also shows the presence of asymmetry
between the amplitudes of the two oscillators along the out-of-phase branch,
just as in the experimental results shown in \figpan{fig:phase_diff}{a}.

As~$K$ increases, the amplitude asymmetry between the oscillators grows
monotonically and the fold of periodics approaches the pitchfork of the
in-phase solution.  At $K \approx 0.14$, the two merge and the pitchfork of
the in-phase solution becomes subcritical. For larger values of $K$, no
stable out-of-phase solutions exist. This suggests that, for sufficiently large~$K$, only the in-phase and anti-phase solutions would be observed in an
experiment.  However, we would still expect bistability between these solutions
due to the presence of an unstable branch of out-of-phase solutions.  Overall,
we find that amplitude asymmetry is strongly associated with bistability of the
phase-locked solutions for non-weak coupling strengths.  This feature cannot be
captured in the phase reduced model \eqref{eq:PhaseOsc} since this approach
disregards amplitude information.

\section{Discussion}
\label{sec:Discussion}

In this article, we investigated transitions between distinct phase-locked
states in a network of two delay-coupled oscillators as the coupling strength
was increased, highlighting the importance of changes in the criticality of said
transitions. One logical question to consider is how these results extend to
networks with more oscillators. Larger networks support a greater variety of
solution types, including partially synchronized cluster
states~\cite{Haugland2021a} and chimera states in which~\cite{Haugland2021} a
portion of the oscillators are phase synchronized, whilst the remaining portion
are not. As such, there is a greater variety of transitions that may occur
between the various states and it would informative to investigate how these how
these change with respect to coupling strength. In our study, we used
DDE-BIFTOOL to analyse the asymptotic solutions of the full system and show how
the criticality of the bifurcations changed. A similar approach could be applied
to study larger networks, however, care must to be taken when discretizing such
systems to ensure that solutions remain accurate but the overall problem remains
numerically tractable.  Moreover, in the case of homogeneous, isotropically
coupled oscillators studied here, the myriad symmetries present in larger
networks can cause numerical difficulties in finding and tracking bifurcation
points.  In this case, additional constraints can be added to the problem
structure to overcome these difficulties~\cite{Krauskopf2023ab}.

A more accurate low-dimensional description of the nonlinear time-delayed system can give more precise insights into the nature of the transition to in-phase or anti-phase synchronized configurations.
The theoretical considerations leading to the results in \fig{fig:Crit}
were based on an ad-hoc phase description with a finite number of harmonics.
Note that for highly nonsinusoidal oscillations---such as relaxation oscillations---a large number of harmonics is required to obtain an accurate description of the phase dynamics, even to first order in coupling strength~\cite{Izhikevich2000,Ashwin2021}.
Computing a phase reduction~\cite{Nakao2015,Pietras2019} explicitly allows to link the phase parameters to the actual physical parameters in the system.
Moreover, higher-order phase reductions remain valid for larger coupling strengths that we would expect in real world experimental systems.
Rigorous reduction approaches for time-delayed systems are only now being developed~\cite{deWolff2024}.
Alternatively, phase-amplitude reduction that include an ``amplitude'' variable in addition to the phase to describe an oscillator's state have proven useful~\cite{Kotani2020}; 
an analysis of phase-amplitude models is beyond the scope of this paper.

The existence of a phase reduction, i.e., that amplitudes are enslaved to phase
variables, is not a contradiction to the asymmetry in amplitudes observed in
\fig{fig:phase_diff} and \fig{fig:bruss}.  While traditional approaches to phase
reduction have focused on deriving approximations for the phase dynamics, a
recent approach based on a parametrization method~\cite{VonderGracht2023a} can
also compute how amplitudes depend on the phases---or, in a more mathematical
language, how an invariant torus is embedded in the state space of the nonlinear
oscillator network.  Hence, this approach can also shed light on the emergent
amplitude dynamics along solutions of the phase equations.

Although our mathematical models do not aim to describe the specific
electrochemical reaction in the experiment, it is still instructive to consider
how well matched the features of the models and the experiment are.  It is
generally impossible in a real-world setting to establish perfectly identical
oscillators, meaning that these systems do not possess the same symmetries as
the mathematical models. However, the discrepancy we here observe is small
($<0.5\%$) and so we consider the oscillator to be approximately identical.  We
also cannot rule out the possibility that we observe in the experiments
long-lived transient behaviour, as opposed to the asymptotic behaviour examined
in our bifurcation analysis. This issue is particularly relevant for the weak
coupling case in which transients may decay over long durations. To mitigate
this, we varied the time delay parameter over a much slower timescale than the
oscillations themselves. In addition, the abrupt transitions we observe for
larger coupling strengths give us confidence that we are sufficiently well
capturing asymptotic dynamics. 
We also note that the findings provide limitations to the extent that the synchronization engineering technique \cite{Kiss2007} can be used for tuning the phase difference between two oscillators \cite{Rusin2010}. In this technique, it is assumed that the feedback gain is sufficiently strong such that the inherent natural frequency difference can be neglected. However, 
in this work, we showed that feedback that is too strong can induce higher order effects that impact the phase dynamics. Therefore, for oscillators with large frequency 
difference, techniques that take advantage of this natural frequency difference, such as phase assignment with resonant entrainment \cite{Zlotnik:2016bja}, are preferable.
Overall, we
expect our results to be relevant to a wide range of applications involving
oscillators networks away from the weak coupling limit.

\section*{Acknowledgements} I.Z.K. acknowledges support from the National Science Foundation (NSF) under Grant No. CHE-1900011.
K.C.A.W gratefully acknowledges the financial support of the EPSRC via grants EP/T017856/1 and EP/V048716/1.

\bibliographystyle{apsrev4-1}

\begin{thebibliography}{41}%
  \makeatletter
  \providecommand \@ifxundefined [1]{%
   \@ifx{#1\undefined}
  }%
  \providecommand \@ifnum [1]{%
   \ifnum #1\expandafter \@firstoftwo
   \else \expandafter \@secondoftwo
   \fi
  }%
  \providecommand \@ifx [1]{%
   \ifx #1\expandafter \@firstoftwo
   \else \expandafter \@secondoftwo
   \fi
  }%
  \providecommand \natexlab [1]{#1}%
  \providecommand \enquote  [1]{``#1''}%
  \providecommand \bibnamefont  [1]{#1}%
  \providecommand \bibfnamefont [1]{#1}%
  \providecommand \citenamefont [1]{#1}%
  \providecommand \href@noop [0]{\@secondoftwo}%
  \providecommand \href [0]{\begingroup \@sanitize@url \@href}%
  \providecommand \@href[1]{\@@startlink{#1}\@@href}%
  \providecommand \@@href[1]{\endgroup#1\@@endlink}%
  \providecommand \@sanitize@url [0]{\catcode `\\12\catcode `\$12\catcode
    `\&12\catcode `\#12\catcode `\^12\catcode `\_12\catcode `\%12\relax}%
  \providecommand \@@startlink[1]{}%
  \providecommand \@@endlink[0]{}%
  \providecommand \url  [0]{\begingroup\@sanitize@url \@url }%
  \providecommand \@url [1]{\endgroup\@href {#1}{\urlprefix }}%
  \providecommand \urlprefix  [0]{URL }%
  \providecommand \Eprint [0]{\href }%
  \providecommand \doibase [0]{http://dx.doi.org/}%
  \providecommand \selectlanguage [0]{\@gobble}%
  \providecommand \bibinfo  [0]{\@secondoftwo}%
  \providecommand \bibfield  [0]{\@secondoftwo}%
  \providecommand \translation [1]{[#1]}%
  \providecommand \BibitemOpen [0]{}%
  \providecommand \bibitemStop [0]{}%
  \providecommand \bibitemNoStop [0]{.\EOS\space}%
  \providecommand \EOS [0]{\spacefactor3000\relax}%
  \providecommand \BibitemShut  [1]{\csname bibitem#1\endcsname}%
  \let\auto@bib@innerbib\@empty
  \bibitem [{\citenamefont {Ashwin}\ \emph {et~al.}(2016)\citenamefont {Ashwin},
    \citenamefont {Coombes},\ and\ \citenamefont {Nicks}}]{Ashwin2016a}%
    \BibitemOpen
    \bibfield  {author} {\bibinfo {author} {\bibfnamefont {P.}~\bibnamefont
    {Ashwin}}, \bibinfo {author} {\bibfnamefont {S.}~\bibnamefont {Coombes}}, \
    and\ \bibinfo {author} {\bibfnamefont {R.}~\bibnamefont {Nicks}},\ }\href
    {\doibase 10.1186/s13408-015-0033-6} {\bibfield  {journal} {\bibinfo
    {journal} {Journal of Mathematical Neuroscience}\ }\textbf {\bibinfo {volume}
    {6}},\ \bibinfo {pages} {1} (\bibinfo {year} {2016})}\BibitemShut {NoStop}%
  \bibitem [{\citenamefont {Bick}\ \emph {et~al.}(2020)\citenamefont {Bick},
    \citenamefont {Goodfellow}, \citenamefont {Laing},\ and\ \citenamefont
    {Martens}}]{Bick2020}%
    \BibitemOpen
    \bibfield  {author} {\bibinfo {author} {\bibfnamefont {C.}~\bibnamefont
    {Bick}}, \bibinfo {author} {\bibfnamefont {M.}~\bibnamefont {Goodfellow}},
    \bibinfo {author} {\bibfnamefont {C.~R.}\ \bibnamefont {Laing}}, \ and\
    \bibinfo {author} {\bibfnamefont {E.~A.}\ \bibnamefont {Martens}},\ }\href
    {\doibase 10.1186/s13408-020-00086-9} {\bibfield  {journal} {\bibinfo
    {journal} {Journal of Mathematical Neuroscience}\ }\textbf {\bibinfo {volume}
    {10}} (\bibinfo {year} {2020}),\ 10.1186/s13408-020-00086-9}\BibitemShut
    {NoStop}%
  \bibitem [{\citenamefont {Filatrella}\ \emph {et~al.}(2008)\citenamefont
    {Filatrella}, \citenamefont {Nielsen},\ and\ \citenamefont
    {Pedersen}}]{Filatrella2008}%
    \BibitemOpen
    \bibfield  {author} {\bibinfo {author} {\bibfnamefont {G.}~\bibnamefont
    {Filatrella}}, \bibinfo {author} {\bibfnamefont {A.~H.}\ \bibnamefont
    {Nielsen}}, \ and\ \bibinfo {author} {\bibfnamefont {N.~F.}\ \bibnamefont
    {Pedersen}},\ }\href {\doibase 10.1140/epjb/e2008-00098-8} {\bibfield
    {journal} {\bibinfo  {journal} {European Physical Journal B}\ }\textbf
    {\bibinfo {volume} {61}},\ \bibinfo {pages} {485} (\bibinfo {year}
    {2008})}\BibitemShut {NoStop}%
  \bibitem [{\citenamefont {D{\"{o}}rfler}\ \emph {et~al.}(2013)\citenamefont
    {D{\"{o}}rfler}, \citenamefont {Chertkov},\ and\ \citenamefont
    {Bullo}}]{Dorfler2013}%
    \BibitemOpen
    \bibfield  {author} {\bibinfo {author} {\bibfnamefont {F.}~\bibnamefont
    {D{\"{o}}rfler}}, \bibinfo {author} {\bibfnamefont {M.}~\bibnamefont
    {Chertkov}}, \ and\ \bibinfo {author} {\bibfnamefont {F.}~\bibnamefont
    {Bullo}},\ }\href {\doibase 10.1073/pnas.1212134110} {\bibfield  {journal}
    {\bibinfo  {journal} {Proceedings of the National Academy of Sciences of the
    United States of America}\ }\textbf {\bibinfo {volume} {110}},\ \bibinfo
    {pages} {2005} (\bibinfo {year} {2013})},\ \Eprint
    {http://arxiv.org/abs/1208.0045} {arXiv:1208.0045} \BibitemShut {NoStop}%
  \bibitem [{\citenamefont {Yan}\ \emph {et~al.}(2007)\citenamefont {Yan},
    \citenamefont {Fu}, \citenamefont {Ren},\ and\ \citenamefont
    {Wang}}]{Yan2007}%
    \BibitemOpen
    \bibfield  {author} {\bibinfo {author} {\bibfnamefont {G.}~\bibnamefont
    {Yan}}, \bibinfo {author} {\bibfnamefont {Z.~Q.}\ \bibnamefont {Fu}},
    \bibinfo {author} {\bibfnamefont {J.}~\bibnamefont {Ren}}, \ and\ \bibinfo
    {author} {\bibfnamefont {W.~X.}\ \bibnamefont {Wang}},\ }\href {\doibase
    10.1103/PhysRevE.75.016108} {\bibfield  {journal} {\bibinfo  {journal}
    {Physical Review E - Statistical, Nonlinear, and Soft Matter Physics}\
    }\textbf {\bibinfo {volume} {75}},\ \bibinfo {pages} {1} (\bibinfo {year}
    {2007})},\ \Eprint {http://arxiv.org/abs/0602137} {arXiv:0602137 [physics]}
    \BibitemShut {NoStop}%
  \bibitem [{\citenamefont {Gross}\ and\ \citenamefont
    {Kevrekidis}(2008)}]{Gross2008}%
    \BibitemOpen
    \bibfield  {author} {\bibinfo {author} {\bibfnamefont {T.}~\bibnamefont
    {Gross}}\ and\ \bibinfo {author} {\bibfnamefont {I.~G.}\ \bibnamefont
    {Kevrekidis}},\ }\href {\doibase 10.1209/0295-5075/82/38004} {\bibfield
    {journal} {\bibinfo  {journal} {Epl}\ }\textbf {\bibinfo {volume} {82}},\
    \bibinfo {pages} {38004} (\bibinfo {year} {2008})},\ \Eprint
    {http://arxiv.org/abs/0702047} {arXiv:0702047 [nlin]} \BibitemShut {NoStop}%
  \bibitem [{\citenamefont {Hoppensteadt}\ and\ \citenamefont
    {Izhikevich}(1997)}]{Hoppensteadt97}%
    \BibitemOpen
    \bibfield  {author} {\bibinfo {author} {\bibfnamefont {F.~C.}\ \bibnamefont
    {Hoppensteadt}}\ and\ \bibinfo {author} {\bibfnamefont {E.~M.}\ \bibnamefont
    {Izhikevich}},\ }\href@noop {} {\emph {\bibinfo {title} {{Weakly Connected
    Neural Networks}}}},\ \bibinfo {series} {Applied Mathematical Sciences}\ No.\
    \bibinfo {number} {126}\ (\bibinfo  {publisher} {Springer-Verlag},\ \bibinfo
    {address} {New York},\ \bibinfo {year} {1997})\BibitemShut {NoStop}%
  \bibitem [{\citenamefont {Zhang}\ \emph {et~al.}(2014)\citenamefont {Zhang},
    \citenamefont {Zou}, \citenamefont {Boccaletti},\ and\ \citenamefont
    {Liu}}]{Zhang2014}%
    \BibitemOpen
    \bibfield  {author} {\bibinfo {author} {\bibfnamefont {X.}~\bibnamefont
    {Zhang}}, \bibinfo {author} {\bibfnamefont {Y.}~\bibnamefont {Zou}}, \bibinfo
    {author} {\bibfnamefont {S.}~\bibnamefont {Boccaletti}}, \ and\ \bibinfo
    {author} {\bibfnamefont {Z.}~\bibnamefont {Liu}},\ }\href {\doibase
    10.1038/srep05200} {\bibfield  {journal} {\bibinfo  {journal} {Scientific
    Reports}\ }\textbf {\bibinfo {volume} {4}},\ \bibinfo {pages} {1} (\bibinfo
    {year} {2014})}\BibitemShut {NoStop}%
  \bibitem [{\citenamefont {Vlasov}\ \emph {et~al.}(2015)\citenamefont {Vlasov},
    \citenamefont {Zou},\ and\ \citenamefont {Pereira}}]{Vlasov2015}%
    \BibitemOpen
    \bibfield  {author} {\bibinfo {author} {\bibfnamefont {V.}~\bibnamefont
    {Vlasov}}, \bibinfo {author} {\bibfnamefont {Y.}~\bibnamefont {Zou}}, \ and\
    \bibinfo {author} {\bibfnamefont {T.}~\bibnamefont {Pereira}},\ }\href
    {\doibase 10.1103/PhysRevE.92.012904} {\bibfield  {journal} {\bibinfo
    {journal} {Physical Review E - Statistical, Nonlinear, and Soft Matter
    Physics}\ }\textbf {\bibinfo {volume} {92}},\ \bibinfo {pages} {1} (\bibinfo
    {year} {2015})},\ \Eprint {http://arxiv.org/abs/1411.6873} {arXiv:1411.6873}
    \BibitemShut {NoStop}%
  \bibitem [{\citenamefont {Nakao}(2016)}]{Nakao2015}%
    \BibitemOpen
    \bibfield  {author} {\bibinfo {author} {\bibfnamefont {H.}~\bibnamefont
    {Nakao}},\ }\href {\doibase 10.1080/00107514.2015.1094987} {\bibfield
    {journal} {\bibinfo  {journal} {Contemporary Physics}\ }\textbf {\bibinfo
    {volume} {57}},\ \bibinfo {pages} {188} (\bibinfo {year} {2016})}\BibitemShut
    {NoStop}%
  \bibitem [{\citenamefont {Pietras}\ and\ \citenamefont
    {Daffertshofer}(2019)}]{Pietras2019}%
    \BibitemOpen
    \bibfield  {author} {\bibinfo {author} {\bibfnamefont {B.}~\bibnamefont
    {Pietras}}\ and\ \bibinfo {author} {\bibfnamefont {A.}~\bibnamefont
    {Daffertshofer}},\ }\href {\doibase 10.1016/j.physrep.2019.06.001} {\bibfield
     {journal} {\bibinfo  {journal} {Phys Rep}\ }\textbf {\bibinfo {volume}
    {819}},\ \bibinfo {pages} {1} (\bibinfo {year} {2019})}\BibitemShut {NoStop}%
  \bibitem [{\citenamefont {Kori}\ \emph {et~al.}(2008)\citenamefont {Kori},
    \citenamefont {Rusin}, \citenamefont {Kiss},\ and\ \citenamefont
    {Hudson}}]{Kori2008}%
    \BibitemOpen
    \bibfield  {author} {\bibinfo {author} {\bibfnamefont {H.}~\bibnamefont
    {Kori}}, \bibinfo {author} {\bibfnamefont {C.~G.}\ \bibnamefont {Rusin}},
    \bibinfo {author} {\bibfnamefont {I.~Z.}\ \bibnamefont {Kiss}}, \ and\
    \bibinfo {author} {\bibfnamefont {J.~L.}\ \bibnamefont {Hudson}},\ }\href
    {\doibase 10.1063/1.2927531} {\bibfield  {journal} {\bibinfo  {journal}
    {Chaos}\ }\textbf {\bibinfo {volume} {18}},\ \bibinfo {pages} {026111}
    (\bibinfo {year} {2008})}\BibitemShut {NoStop}%
  \bibitem [{\citenamefont {Kiss}(2018)}]{Kiss2018}%
    \BibitemOpen
    \bibfield  {author} {\bibinfo {author} {\bibfnamefont {I.~Z.}\ \bibnamefont
    {Kiss}},\ }\href {\doibase 10.1016/j.coche.2018.02.006} {\bibfield  {journal}
    {\bibinfo  {journal} {Current Opinion in Chemical Engineering}\ }\textbf
    {\bibinfo {volume} {21}},\ \bibinfo {pages} {1} (\bibinfo {year}
    {2018})}\BibitemShut {NoStop}%
  \bibitem [{\citenamefont {Călugăru}\ \emph {et~al.}(2020)\citenamefont
    {Călugăru}, \citenamefont {Totz}, \citenamefont {Martens},\ and\
    \citenamefont {Engel}}]{Calugaru2018}%
    \BibitemOpen
    \bibfield  {author} {\bibinfo {author} {\bibfnamefont {D.}~\bibnamefont
    {Călugăru}}, \bibinfo {author} {\bibfnamefont {J.~F.}\ \bibnamefont
    {Totz}}, \bibinfo {author} {\bibfnamefont {E.~A.}\ \bibnamefont {Martens}}, \
    and\ \bibinfo {author} {\bibfnamefont {H.}~\bibnamefont {Engel}},\ }\href
    {\doibase 10.1126/sciadv.abb2637} {\bibfield  {journal} {\bibinfo  {journal}
    {Science Advances}\ }\textbf {\bibinfo {volume} {6}},\ \bibinfo {pages}
    {eabb2637} (\bibinfo {year} {2020})}\BibitemShut {NoStop}%
  \bibitem [{\citenamefont {B{\"{o}}rgers}(2023)}]{Borgers2023}%
    \BibitemOpen
    \bibfield  {author} {\bibinfo {author} {\bibfnamefont {C.}~\bibnamefont
    {B{\"{o}}rgers}},\ }\href {\doibase 10.1016/j.exco.2023.100120} {\bibfield
    {journal} {\bibinfo  {journal} {Examples and Counterexamples}\ }\textbf
    {\bibinfo {volume} {4}},\ \bibinfo {pages} {100120} (\bibinfo {year}
    {2023})}\BibitemShut {NoStop}%
  \bibitem [{\citenamefont {Guckenheimer}\ and\ \citenamefont
    {Worfolk}(1992)}]{Guckenheimer1992}%
    \BibitemOpen
    \bibfield  {author} {\bibinfo {author} {\bibfnamefont {J.}~\bibnamefont
    {Guckenheimer}}\ and\ \bibinfo {author} {\bibfnamefont {P.}~\bibnamefont
    {Worfolk}},\ }\href {\doibase 10.1088/0951-7715/5/6/001} {\bibfield
    {journal} {\bibinfo  {journal} {Nonlinearity}\ }\textbf {\bibinfo {volume}
    {5}},\ \bibinfo {pages} {1211} (\bibinfo {year} {1992})}\BibitemShut
    {NoStop}%
  \bibitem [{\citenamefont {Skardal}\ and\ \citenamefont
    {Arenas}(2020)}]{Skardal2020}%
    \BibitemOpen
    \bibfield  {author} {\bibinfo {author} {\bibfnamefont {P.~S.}\ \bibnamefont
    {Skardal}}\ and\ \bibinfo {author} {\bibfnamefont {A.}~\bibnamefont
    {Arenas}},\ }\href {\doibase 10.1038/s42005-020-00485-0} {\bibfield
    {journal} {\bibinfo  {journal} {Communications Physics}\ }\textbf {\bibinfo
    {volume} {3}} (\bibinfo {year} {2020}),\ 10.1038/s42005-020-00485-0},\
    \Eprint {http://arxiv.org/abs/1909.08057} {arXiv:1909.08057} \BibitemShut
    {NoStop}%
  \bibitem [{\citenamefont {Bick}\ \emph {et~al.}(2023)\citenamefont {Bick},
    \citenamefont {Gross}, \citenamefont {Harrington},\ and\ \citenamefont
    {Schaub}}]{Bick2023}%
    \BibitemOpen
    \bibfield  {author} {\bibinfo {author} {\bibfnamefont {C.}~\bibnamefont
    {Bick}}, \bibinfo {author} {\bibfnamefont {E.}~\bibnamefont {Gross}},
    \bibinfo {author} {\bibfnamefont {H.~A.}\ \bibnamefont {Harrington}}, \ and\
    \bibinfo {author} {\bibfnamefont {M.~T.}\ \bibnamefont {Schaub}},\
    }\href@noop {} {\bibfield  {journal} {\bibinfo  {journal} {SIAM Review}\
    }\textbf {\bibinfo {volume} {65}},\ \bibinfo {pages} {686} (\bibinfo {year}
    {2023})}\BibitemShut {NoStop}%
  \bibitem [{\citenamefont {Wedgwood}\ \emph {et~al.}(2013)\citenamefont
    {Wedgwood}, \citenamefont {Lin}, \citenamefont {Thul},\ and\ \citenamefont
    {Coombes}}]{Wedgwood2013}%
    \BibitemOpen
    \bibfield  {author} {\bibinfo {author} {\bibfnamefont {K.~C.~A.}\
    \bibnamefont {Wedgwood}}, \bibinfo {author} {\bibfnamefont {K.~K.}\
    \bibnamefont {Lin}}, \bibinfo {author} {\bibfnamefont {R.}~\bibnamefont
    {Thul}}, \ and\ \bibinfo {author} {\bibfnamefont {S.}~\bibnamefont
    {Coombes}},\ }\href {\doibase 10.1186/2190-8567-3-2} {\bibfield  {journal}
    {\bibinfo  {journal} {Journal of Mathematical Neuroscience}\ }\textbf
    {\bibinfo {volume} {3}},\ \bibinfo {pages} {2} (\bibinfo {year}
    {2013})}\BibitemShut {NoStop}%
  \bibitem [{\citenamefont {Wilson}\ and\ \citenamefont
    {Moehlis}(2016)}]{Wilson2016}%
    \BibitemOpen
    \bibfield  {author} {\bibinfo {author} {\bibfnamefont {D.}~\bibnamefont
    {Wilson}}\ and\ \bibinfo {author} {\bibfnamefont {J.}~\bibnamefont
    {Moehlis}},\ }\href {\doibase 10.1103/PhysRevE.94.052213} {\bibfield
    {journal} {\bibinfo  {journal} {Physical Review E}\ }\textbf {\bibinfo
    {volume} {94}},\ \bibinfo {pages} {052213} (\bibinfo {year}
    {2016})}\BibitemShut {NoStop}%
  \bibitem [{\citenamefont {Letson}\ and\ \citenamefont
    {Rubin}(2018)}]{Letson2018}%
    \BibitemOpen
    \bibfield  {author} {\bibinfo {author} {\bibfnamefont {B.}~\bibnamefont
    {Letson}}\ and\ \bibinfo {author} {\bibfnamefont {J.~E.}\ \bibnamefont
    {Rubin}},\ }\href {\doibase 10.1137/18m1186617} {\bibfield  {journal}
    {\bibinfo  {journal} {SIAM Journal on Applied Dynamical Systems}\ }\textbf
    {\bibinfo {volume} {17}},\ \bibinfo {pages} {2414} (\bibinfo {year}
    {2018})}\BibitemShut {NoStop}%
  \bibitem [{\citenamefont {Wilson}\ and\ \citenamefont
    {Ermentrout}(2019)}]{Wilson2019}%
    \BibitemOpen
    \bibfield  {author} {\bibinfo {author} {\bibfnamefont {D.}~\bibnamefont
    {Wilson}}\ and\ \bibinfo {author} {\bibfnamefont {B.}~\bibnamefont
    {Ermentrout}},\ }\href {\doibase 10.1103/PhysRevLett.123.164101} {\bibfield
    {journal} {\bibinfo  {journal} {Physical Review Letters}\ }\textbf {\bibinfo
    {volume} {123}},\ \bibinfo {pages} {164101} (\bibinfo {year}
    {2019})}\BibitemShut {NoStop}%
  \bibitem [{\citenamefont {Kotani}\ \emph {et~al.}(2012)\citenamefont {Kotani},
    \citenamefont {Yamaguchi}, \citenamefont {Ogawa}, \citenamefont {Jimbo},
    \citenamefont {Nakao},\ and\ \citenamefont {Ermentrout}}]{Kotani2012}%
    \BibitemOpen
    \bibfield  {author} {\bibinfo {author} {\bibfnamefont {K.}~\bibnamefont
    {Kotani}}, \bibinfo {author} {\bibfnamefont {I.}~\bibnamefont {Yamaguchi}},
    \bibinfo {author} {\bibfnamefont {Y.}~\bibnamefont {Ogawa}}, \bibinfo
    {author} {\bibfnamefont {Y.}~\bibnamefont {Jimbo}}, \bibinfo {author}
    {\bibfnamefont {H.}~\bibnamefont {Nakao}}, \ and\ \bibinfo {author}
    {\bibfnamefont {G.~B.}\ \bibnamefont {Ermentrout}},\ }\href {\doibase
    10.1103/PhysRevLett.109.044101} {\bibfield  {journal} {\bibinfo  {journal}
    {Physical Review Letters}\ }\textbf {\bibinfo {volume} {109}},\ \bibinfo
    {pages} {044101} (\bibinfo {year} {2012})}\BibitemShut {NoStop}%
  \bibitem [{\citenamefont {Wilson}(2020)}]{Wilson2020b}%
    \BibitemOpen
    \bibfield  {author} {\bibinfo {author} {\bibfnamefont {D.}~\bibnamefont
    {Wilson}},\ }\href {\doibase 10.1063/1.5126122} {\bibfield  {journal}
    {\bibinfo  {journal} {Chaos}\ }\textbf {\bibinfo {volume} {30}} (\bibinfo
    {year} {2020}),\ 10.1063/1.5126122}\BibitemShut {NoStop}%
  \bibitem [{\citenamefont {Sakaguchi}\ and\ \citenamefont
    {Kuramoto}(1986)}]{Sakaguchi1986}%
    \BibitemOpen
    \bibfield  {author} {\bibinfo {author} {\bibfnamefont {H.}~\bibnamefont
    {Sakaguchi}}\ and\ \bibinfo {author} {\bibfnamefont {Y.}~\bibnamefont
    {Kuramoto}},\ }\href {\doibase 10.1143/PTP.76.576} {\bibfield  {journal}
    {\bibinfo  {journal} {Progress of Theoretical Physics}\ }\textbf {\bibinfo
    {volume} {76}},\ \bibinfo {pages} {576} (\bibinfo {year} {1986})}\BibitemShut
    {NoStop}%
  \bibitem [{\citenamefont {Rusin}\ \emph {et~al.}(2010)\citenamefont {Rusin},
    \citenamefont {Kori}, \citenamefont {Kiss},\ and\ \citenamefont
    {Hudson}}]{Rusin2010}%
    \BibitemOpen
    \bibfield  {author} {\bibinfo {author} {\bibfnamefont {C.~G.}\ \bibnamefont
    {Rusin}}, \bibinfo {author} {\bibfnamefont {H.}~\bibnamefont {Kori}},
    \bibinfo {author} {\bibfnamefont {I.~Z.}\ \bibnamefont {Kiss}}, \ and\
    \bibinfo {author} {\bibfnamefont {J.~L.}\ \bibnamefont {Hudson}},\ }\href
    {\doibase 10.1098/rsta.2010.0032} {\bibfield  {journal} {\bibinfo  {journal}
    {Philos T Ro Soc A}\ }\textbf {\bibinfo {volume} {368}},\ \bibinfo {pages}
    {2189} (\bibinfo {year} {2010})}\BibitemShut {NoStop}%
  \bibitem [{\citenamefont {Kuehn}\ and\ \citenamefont {Bick}(2021)}]{Kuehn2020}%
    \BibitemOpen
    \bibfield  {author} {\bibinfo {author} {\bibfnamefont {C.}~\bibnamefont
    {Kuehn}}\ and\ \bibinfo {author} {\bibfnamefont {C.}~\bibnamefont {Bick}},\
    }\href {\doibase 10.1126/sciadv.abe3824} {\bibfield  {journal} {\bibinfo
    {journal} {Science Advances}\ }\textbf {\bibinfo {volume} {7}},\ \bibinfo
    {pages} {eabe3824} (\bibinfo {year} {2021})}\BibitemShut {NoStop}%
  \bibitem [{Note1()}]{Note1}%
    \BibitemOpen
    \bibinfo {note} {For each parameter~$\alpha _n$, we solved~\protect \textup
    {\hbox {\mathsurround \z@ \protect \normalfont (\ignorespaces \ref
    {eq:PhaseDiff}\unskip \@@italiccorr )}}, numerically for $T=5,000$ time units
    with initial condition being~$\psi (T)$ for parameter $\alpha _{n\pm 1}$ plus
    a small random perturbation.}\BibitemShut {Stop}%
  \bibitem [{\citenamefont {Kiss}\ \emph {et~al.}(2006)\citenamefont {Kiss},
    \citenamefont {Kazsu},\ and\ \citenamefont
    {G{\'a}sp{\'a}r}}]{kiss_tracking_2006}%
    \BibitemOpen
    \bibfield  {author} {\bibinfo {author} {\bibfnamefont {I.~Z.}\ \bibnamefont
    {Kiss}}, \bibinfo {author} {\bibfnamefont {Z.}~\bibnamefont {Kazsu}}, \ and\
    \bibinfo {author} {\bibfnamefont {V.}~\bibnamefont {G{\'a}sp{\'a}r}},\ }\href
    {\doibase 10.1063/1.2219702} {\bibfield  {journal} {\bibinfo  {journal}
    {Chaos}\ }\textbf {\bibinfo {volume} {16}},\ \bibinfo {pages} {033109}
    (\bibinfo {year} {2006})}\BibitemShut {NoStop}%
  \bibitem [{\citenamefont {Bick}\ \emph {et~al.}(2017)\citenamefont {Bick},
    \citenamefont {Sebek},\ and\ \citenamefont {Kiss}}]{bick2017robust}%
    \BibitemOpen
    \bibfield  {author} {\bibinfo {author} {\bibfnamefont {C.}~\bibnamefont
    {Bick}}, \bibinfo {author} {\bibfnamefont {M.}~\bibnamefont {Sebek}}, \ and\
    \bibinfo {author} {\bibfnamefont {I.~Z.}\ \bibnamefont {Kiss}},\ }\href@noop
    {} {\bibfield  {journal} {\bibinfo  {journal} {Phys. Rev. Lett.}\ }\textbf
    {\bibinfo {volume} {119}},\ \bibinfo {pages} {168301} (\bibinfo {year}
    {2017})}\BibitemShut {NoStop}%
  \bibitem [{\citenamefont {Hansel}\ \emph {et~al.}(1993)\citenamefont {Hansel},
    \citenamefont {Mato},\ and\ \citenamefont {Meunier}}]{Hansel1993}%
    \BibitemOpen
    \bibfield  {author} {\bibinfo {author} {\bibfnamefont {D.}~\bibnamefont
    {Hansel}}, \bibinfo {author} {\bibfnamefont {G.}~\bibnamefont {Mato}}, \ and\
    \bibinfo {author} {\bibfnamefont {C.}~\bibnamefont {Meunier}},\ }\href
    {\doibase 10.1103/PhysRevE.48.3470} {\bibfield  {journal} {\bibinfo
    {journal} {Physical Review E}\ }\textbf {\bibinfo {volume} {48}},\ \bibinfo
    {pages} {3470} (\bibinfo {year} {1993})}\BibitemShut {NoStop}%
  \bibitem [{\citenamefont {Haugland}\ \emph {et~al.}(2021)\citenamefont
    {Haugland}, \citenamefont {Tosolini},\ and\ \citenamefont
    {Krischer}}]{Haugland2021a}%
    \BibitemOpen
    \bibfield  {author} {\bibinfo {author} {\bibfnamefont {S.~W.}\ \bibnamefont
    {Haugland}}, \bibinfo {author} {\bibfnamefont {A.}~\bibnamefont {Tosolini}},
    \ and\ \bibinfo {author} {\bibfnamefont {K.}~\bibnamefont {Krischer}},\
    }\href {\doibase 10.1038/s41467-021-25907-7} {\bibfield  {journal} {\bibinfo
    {journal} {Nature Communications}\ }\textbf {\bibinfo {volume} {12}}
    (\bibinfo {year} {2021}),\ 10.1038/s41467-021-25907-7}\BibitemShut {NoStop}%
  \bibitem [{\citenamefont {Haugland}(2021)}]{Haugland2021}%
    \BibitemOpen
    \bibfield  {author} {\bibinfo {author} {\bibfnamefont {S.~W.}\ \bibnamefont
    {Haugland}},\ }\href {\doibase 10.1088/2632-072X/ac0810} {\bibfield  {journal} {\bibinfo
    {journal} {Journal of Physics: Complexity}\ }\textbf {\bibinfo {volume} {12}}, {\bibinfo {number} {3}} (\bibinfo {year} {2021}),\ 10.1038/s41467-021-25907-7}
    \BibitemShut {NoStop}%
  \bibitem [{\citenamefont {Krauskopf}\ and\ \citenamefont
    {Sieber}(2023)}]{Krauskopf2023ab}%
    \BibitemOpen
    \bibfield  {author} {\bibinfo {author} {\bibfnamefont {B.}~\bibnamefont
    {Krauskopf}}\ and\ \bibinfo {author} {\bibfnamefont {J.}~\bibnamefont
    {Sieber}},\ }\enquote {\bibinfo {title} {Bifurcation analysis of systems with
    delays: Methods and their use in applications},}\ \ (\bibinfo  {publisher}
    {Springer Cham},\ \bibinfo {year} {2023})\ pp.\ \bibinfo {pages}
    {195--245}\BibitemShut {NoStop}%
  \bibitem [{\citenamefont {Izhikevich}(2000)}]{Izhikevich2000}%
    \BibitemOpen
    \bibfield  {author} {\bibinfo {author} {\bibfnamefont {E.~M.}\ \bibnamefont
    {Izhikevich}},\ }\href {\doibase 10.1137/S0036139999351001} {\bibfield
    {journal} {\bibinfo  {journal} {SIAM Journal on Applied Mathematics}\
    }\textbf {\bibinfo {volume} {60}},\ \bibinfo {pages} {1789} (\bibinfo {year}
    {2000})}\BibitemShut {NoStop}%
  \bibitem [{\citenamefont {Ashwin}\ \emph {et~al.}(2021)\citenamefont {Ashwin},
    \citenamefont {Bick},\ and\ \citenamefont {Poignard}}]{Ashwin2021}%
    \BibitemOpen
    \bibfield  {author} {\bibinfo {author} {\bibfnamefont {P.}~\bibnamefont
    {Ashwin}}, \bibinfo {author} {\bibfnamefont {C.}~\bibnamefont {Bick}}, \ and\
    \bibinfo {author} {\bibfnamefont {C.}~\bibnamefont {Poignard}},\ }\href
    {\doibase 10.1063/5.0063423} {\bibfield  {journal} {\bibinfo  {journal}
    {Chaos}\ }\textbf {\bibinfo {volume} {31}},\ \bibinfo {pages} {093132}
    (\bibinfo {year} {2021})},\ \Eprint {http://arxiv.org/abs/2107.07152}
    {2107.07152} \BibitemShut {NoStop}%
  \bibitem [{\citenamefont {Bick}\ \emph {et~al.}(2024)\citenamefont {Bick},
    \citenamefont {Rink},\ and\ \citenamefont {de~Wolff}}]{deWolff2024}%
    \BibitemOpen
    \bibfield  {author} {\bibinfo {author} {\bibfnamefont {C.}~\bibnamefont
    {Bick}}, \bibinfo {author} {\bibfnamefont {B.}~\bibnamefont {Rink}}, \ and\
    \bibinfo {author} {\bibfnamefont {B.}~\bibnamefont {de~Wolff}},\ }\href@noop
    {} {\enquote {\bibinfo {title} {{Phase reductions for delay-coupled
    oscillators}},}\ } (\bibinfo {year} {2024}),\ \bibinfo {note} {{In
    Prep}}\BibitemShut {NoStop}%
  \bibitem [{\citenamefont {Kotani}\ \emph {et~al.}(2020)\citenamefont {Kotani},
    \citenamefont {Ogawa}, \citenamefont {Shirasaka}, \citenamefont {Akao},
    \citenamefont {Jimbo},\ and\ \citenamefont {Nakao}}]{Kotani2020}%
    \BibitemOpen
    \bibfield  {author} {\bibinfo {author} {\bibfnamefont {K.}~\bibnamefont
    {Kotani}}, \bibinfo {author} {\bibfnamefont {Y.}~\bibnamefont {Ogawa}},
    \bibinfo {author} {\bibfnamefont {S.}~\bibnamefont {Shirasaka}}, \bibinfo
    {author} {\bibfnamefont {A.}~\bibnamefont {Akao}}, \bibinfo {author}
    {\bibfnamefont {Y.}~\bibnamefont {Jimbo}}, \ and\ \bibinfo {author}
    {\bibfnamefont {H.}~\bibnamefont {Nakao}},\ }\href {\doibase
    10.1103/PhysRevResearch.2.033106} {\bibfield  {journal} {\bibinfo  {journal}
    {Physical Review Research}\ }\textbf {\bibinfo {volume} {2}},\ \bibinfo
    {pages} {033106} (\bibinfo {year} {2020})}\BibitemShut {NoStop}%
  \bibitem [{\citenamefont {von~der Gracht}\ \emph {et~al.}(2023)\citenamefont
    {von~der Gracht}, \citenamefont {Nijholt},\ and\ \citenamefont
    {Rink}}]{VonderGracht2023a}%
    \BibitemOpen
    \bibfield  {author} {\bibinfo {author} {\bibfnamefont {S.}~\bibnamefont
    {von~der Gracht}}, \bibinfo {author} {\bibfnamefont {E.}~\bibnamefont
    {Nijholt}}, \ and\ \bibinfo {author} {\bibfnamefont {B.}~\bibnamefont
    {Rink}},\ }\href {http://arxiv.org/abs/2306.03320} {\bibfield  {journal}
    {\bibinfo  {journal} {arXiv:2306.03320}\ } (\bibinfo {year}
    {2023})}\BibitemShut {NoStop}%
  \bibitem [{\citenamefont {Kiss}\ \emph {et~al.}(2007)\citenamefont {Kiss},
    \citenamefont {Rusin}, \citenamefont {Kori},\ and\ \citenamefont
    {Hudson}}]{Kiss2007}%
    \BibitemOpen
    \bibfield  {author} {\bibinfo {author} {\bibfnamefont {I.~Z.}\ \bibnamefont
    {Kiss}}, \bibinfo {author} {\bibfnamefont {C.~G.}\ \bibnamefont {Rusin}},
    \bibinfo {author} {\bibfnamefont {H.}~\bibnamefont {Kori}}, \ and\ \bibinfo
    {author} {\bibfnamefont {J.~L.}\ \bibnamefont {Hudson}},\ }\href {\doibase
    10.1126/science.1140858} {\bibfield  {journal} {\bibinfo  {journal}
    {Science}\ }\textbf {\bibinfo {volume} {316}},\ \bibinfo {pages} {1886}
    (\bibinfo {year} {2007})}\BibitemShut {NoStop}%
  \bibitem [{\citenamefont {Zlotnik}\ \emph {et~al.}(2016)\citenamefont
    {Zlotnik}, \citenamefont {Nagao}, \citenamefont {Kiss},\ and\ \citenamefont
    {Li}}]{Zlotnik:2016bja}%
    \BibitemOpen
    \bibfield  {author} {\bibinfo {author} {\bibfnamefont {A.}~\bibnamefont
    {Zlotnik}}, \bibinfo {author} {\bibfnamefont {R.}~\bibnamefont {Nagao}},
    \bibinfo {author} {\bibfnamefont {I.~Z.}\ \bibnamefont {Kiss}}, \ and\
    \bibinfo {author} {\bibfnamefont {J.-S.}\ \bibnamefont {Li}},\ }\href
    {\doibase 10.1038/ncomms10788} {\bibfield  {journal} {\bibinfo  {journal}
    {Nature Communications}\ }\textbf {\bibinfo {volume} {7}},\ \bibinfo {pages}
    {10788} (\bibinfo {year} {2016})}\BibitemShut {NoStop}%
  \end{thebibliography}

\appendix

\section{Pipeline for numerical bifurcation analysis}

The continuation of solutions with respect to a single, common delay performed
in \fig{fig:bruss} requires the use of DDE-BIFTOOL's functionality to include
constraints between parameters. Practically speaking, we apply a
\textit{sys\_cond} function that fixes the difference $\tau_1 - \tau_2$ to be
constant and absorb this into the common delay $\tau$.  Following this, we use
the following pipeline:

\begin{enumerate}
  \item Compute the isolated periodic orbit for \eqref{eq:bruss}
    when $K=0$.
  \item Set $K = \overline{K} = \Delta K = 0.01$.
  \item Set an initial condition for the coupled system by shifting the
    periodic orbit of one of the oscillators by $\alpha \in [0,\pi]$. In
    practice, this is done by representing the orbit by its Fourier series and
    then multiplying each of the Fourier coefficients by
    $\ee^{-i\alpha/T}$.
  \item Converge an initial periodic solution for the coupled system using a bespoke
    Matlab function. Repeat for the in-phase branch ($\alpha = 0$), the
    anti-phase branch ($\alpha = \pi$), and the out-of-phase branch ($\alpha \in
    (0,\pi))$.
  \item Continue each branch of periodic solutions over the range $\tau
    \in [0, T/2]$.
  \item Select a point along each branch and continue said branch over the range
    $K \in [\overline{K}, \overline{K} + \Delta K]$.
  \item Increment $\overline{K}$ by $\Delta K$ and repeat steps 5-7 until the maximum value of $K$
    is reached.
  \item Finally, extract summary statistics along each branch, including:
    \begin{itemize}
      \item Linear stability as determined by Floquet multipliers;
      \item Phase difference, $\psi$, assessed by computing the absolute time difference between
        the peaks in $x_i$ for the two oscillators and scaling this by
        $2\pi/\widetilde{T}$
        where $\widetilde{T}$ is the period of the periodic solution;
      \item The relative amplitude asymmetry between the orbits of the two electrodes,
        $\Delta A = |A_1 - A_2|/\min\{A_1,A_2\}$, where $A_i = x_i^\text{max} -
        x_i^\text{min}$.
    \end{itemize}
\end{enumerate}

\end{document}